\documentclass[11pt]{article}
    
\usepackage{amsmath,amsfonts,amsthm,amscd,amssymb,graphicx}
\usepackage{mathabx}

\usepackage[utf8]{inputenc}

\usepackage{subfigure}
\usepackage{xcolor}

\numberwithin{equation}{section}

\usepackage{hyperref}

\newtheorem{theorem}{Theorem}[section]

\newtheorem{lemma}[theorem]{Lemma}

\newtheorem{proposition}[theorem]{Proposition}

\newtheorem{corollary}[theorem]{Corollary}

\newcommand{\RR}{\mathbb{R}}

\newcommand{\cL}{\mathcal{L}}

\newcommand{\cE}{\mathcal{E}}
\newcommand{\cB}{\mathcal{B}}

\newcommand{\TT}{\widetilde{T}}

\newcommand{\DD}{\mathbb{D}}
\newcommand{\PP}{\mathbb{P}}

\begin{document}

\title{Propagation of chaos with singular potentials: A weighted duality approach} 

\author{
Phan Th\`anh Nam\footnotemark[1] \and Toan T. Nguyen\footnotemark[2]
}

\footnotetext[1]{LMU Munich, Department of Mathematics, Theresienstr. 39, 80333 Munich, Germany. Email: nam@math.lmu.de.}

\footnotetext[2]{Penn State University, Department of Mathematics, State College, PA 16802, USA. Email: nguyen@math.psu.edu.}

\maketitle

\begin{abstract}
We establish the propagation of chaos and derive the Vlasov equation from a many-body $N$-particle system in the mean-field regime with singular interaction forces bounded by $|x|^{-\alpha-1}$ for $0\le \alpha <(d-1)/2$ in dimensions $d\ge 2$, which in particular include the full range of sub-Coulomb interactions in three dimensions as well as the Coulomb interaction in two dimensions. Our proof is based on an energy-based framework for deriving weighted uniform bounds on cumulants that builds upon the recent duality method developed by Bresch--Duerinckx--Jabin in combination with a new Morawetz-type estimate, allowing us to remove the restriction to locally square-integrable interaction force kernels and handle more singular regimes.
\end{abstract}


\tableofcontents


\section{Introduction}


The rigorous derivation of effective macroscopic equations from the microscopic dynamics of interacting many-body systems stands as a fundamental problem in mathematical physics. In this direction, the derivation of the Vlasov equation starting directly from the deterministic Newton equations governing an $N$-particle system is a great challenge for the Coulomb interaction in the kinetic theory of plasmas. Heuristically, the formal derivation relies on the assumption that the constituent particles are statistically independent, naturally leading to the mean-field approximation where the total interaction force experienced by a single particle can be effectively replaced by a self-consistent macroscopic field generated by the system's one-body reduced spatial density. The mathematical difficulty in justifying this limit lies in proving rigorously that this statistical independence actually holds over macroscopic time scales, a property canonically referred to as the propagation of chaos.

The study of mean-field limits is a classical subject, starting from the pioneering works \cite{McKean, BH, Dob} that treat smooth interaction kernels; see also the excellent surveys \cite{Golse, Jabinsurvey} on the subject. The problem has recently been studied intensively, focusing on singular interaction kernels, especially the Riesz type potentials \eqref{Riesz}. For the first order particle systems, we mention \cite{Serfaty} which treats the full range of  Riesz kernels, see also \cite{JabinWang, Guillin, Rosenzweig, Nguyen,HS-26}. On the other hand, fewer works on the mean-field limit for the second order systems are available, however we mention the recent important works \cite{BJS, Jabin,Jabin2} that treat a class of singular interaction kernels. We also refer to \cite{DP-17} for the derivation of Vlasov-Poisson equation from potentials with suitable cut-off, and \cite{HP-25,HP-25b} for recent refinements.

Our work is inspired by the recent development \cite{Jabin} of Bresch, Duerinckx and Jabin, who introduced a duality method to establish mean-field limits by working with a dual BBGKY hierarchy for cumulants.  More precisely, by establishing $L^2$-uniform bounds on the cumulants of the many-body dynamics, they are able to show that particle independence, if present at the initial time, is dynamically preserved despite the continuous accumulation of microscopic correlations induced by the interactions. Consequently, they rigorously derive the mean-field limit for all locally square-integrable interaction force kernels.

In the present paper, we propose an energy-based framework which combines the duality argument in \cite{Jabin} with a new weighted uniform bound on cumulants via a Morawetz-type estimate, allowing us to handle more singular regimes. Let us now precise the mathematical setting of the microscopic dynamics as well as the effective macroscopic equations, and then state precisely our result.

\subsection{Newton and Vlasov equations}

We are interested in the dynamics of $N$ interacting particles $\{(x_i(t), v_i(t))\}_{i=1}^N$ in the phase space $\DD:=\RR^d\times \RR^d$, with $d\ge 2$, following  the classical Newtonian dynamics subject to Brownian forces in the mean-field regime:
\begin{equation}\label{Nsys}
\left\{ 
\begin{aligned}
d x_i &= v_i dt
\\
dv_i &= \frac{1}{N-1}\sum^N_{j \not =i} K(x_i - x_j) dt + \sqrt{2\nu_N } dW_i
\end{aligned}
\right.
\end{equation}
for each $i = 1,\dots, N$, where $\{W_i\}$ are  independent Brownian motions with a temperature $0<\nu_N\ll1$ which is vanishing in the limit $N\to \infty$. The particles interact  through the interaction force $K: \mathbb{R}^d \to \mathbb{R}^d$ satisfying 
\begin{equation}\label{Riesz}
K(-x) = - K(x),\qquad |K(x)| \le  \frac{C_K}{|x|^{\alpha+1}}
\end{equation}
with constants $C_K>0$ and $0\le \alpha < \frac{d-1}{2}$. The reader may keep in mind the form $K(x)=-\nabla w(x)$ with an even potential $w:\mathbb{R}^d\to \mathbb{R}$, and the condition  \eqref{Riesz} includes all the sub-Coulomb interaction potentials $w(x)=|x|^{-\alpha}$ in three dimensions as well as the Coulomb potential $w(x)=-\log |x|$ in two dimensions.  The sign of the interaction plays no role in our analysis. 

Our aim is to derive the effective equation describing the macroscopic behavior of the $N$-body dynamics  $\{(x_i(t), v_i(t))\}_{i=1}^N$ when $N\to \infty$. To this end, we consider the equivalent formulation of \eqref{Nsys} where the dynamics of $N$-particle density distribution function $F_N(t,z_{[N]})$, with $z_{[N]} = (z_1,\cdots, z_N)$ and $z_i= (x_i,v_i)$, solving the Liouville equation 
\begin{equation}\label{NVlasov}
 \partial_t F_N + \cL_N F_N = 0
 \end{equation} 
in $N$-particle phase space $\DD^N$, with 
\begin{equation}\label{def-LN}
\cL_N = \sum_{i=1}^N v_i \cdot \nabla_{x_i} +  \frac{1}{N-1}\sum^N_{i\not =j} K(x_i - x_j)\cdot \nabla_{v_i} -\nu_N  \sum_i^N \Delta_{v_i} 
\end{equation}
acting on functions of $z_{[N]} = (z_1,\cdots, z_N)$ in $\DD^N$. We shall focus on symmetric functions $F_N(t,z_{[N]})$, namely \begin{equation}\label{symmF}F_N(t,z_{[N]}) = F_N(t,z_{\sigma_N})\end{equation}
for any permutation $\sigma_N$ of $[N] = \{1,\cdots,N\}$. 

We shall consider the macroscopic behavior of the dynamics of \eqref{NVlasov} under the assumption that the initial data is a factorized state, namely
\begin{equation}
\label{Ndata}F_N(0,z_{[N]}) = f_0^{\otimes N}(z_{[N]}) = \prod_{j=1}^N f_0(z_j)
\end{equation}
for a smooth initial profile $f_0(x,v)$. We expect that since the particles in the initial state are independent, this property remains approximately true for the state $F_N(t,z_{[N]})$ at any finite time $t>0$, namely $F_N(t,z_{[N]})\approx f(t)^{\otimes N}$ in an appropriate sense. By formally inserting this factorization ansatz into \eqref{NVlasov} and applying the law of large numbers under the independent distribution governed by $f$, the empirical force experienced by the $i$-th particle can be replaced by its mean-field approximation: namely, 
$$
\frac{1}{N-1}\sum_{j \neq i}^N K(x_i - x_j) \approx (K \star_x \rho_f) (x_i),\quad \rho_f = \int_{\mathbb{R}^d} f \;dv. 
$$
Consequently, we arrive at the Vlasov equation
\begin{equation}\label{VP}
(\partial_t + \cL) f = 0,\qquad \cL = v\cdot \nabla_x + K[f] \cdot \nabla_v,\quad K[f] = K\star_x \rho_f
\end{equation}
with the same initial data $f(0,x,v) = f_0(x,v)$.

Note that the approximation $F_N(t,z_{[N]})\approx f(t)^{\otimes N}$ is highly sensitive to microscopic correlations, and hence we can only expect to establish such an approximation in a collective sense. To formulate this mathematically, we introduce the $k^{th}$-marginals defined by
\begin{equation}\label{def-FNk}
F_{N,k}(t, z_{[k]}) = \int_{\DD^{N-k}} F_N(t, z_{[N]}) \; dz_{k+1}\cdots dz_N,
\end{equation}
for $1 \le k < N$. In the mean-field limit as $N\to \infty$, we expect that $F_{N,k}(t) \approx f(t)^{\otimes k}$ for any fixed $k\ge 1$, which provides the precise meaning of the approximate factorization (or molecular chaos).

\subsection{Main results}

For simplicity, we shall assume that the initial data $f_0(x,v)\ge 0$ is smooth and compactly supported in phase space. Under this assumption, it is well known that in dimensions $d\in \{2,3\}$, the Vlasov equation \eqref{VP}, including the repulsive Coulomb interaction, admits global smooth solutions that remain compactly supported for all time \cite{LP91}. In higher dimensions, such a smooth, compactly supported solution exists at least up to a finite time. 

Moreover, note that for every parameter $\epsilon \in (0,1)$, the $\epsilon$-power function $h=f(t,z)^\epsilon$ satisfies the same linear transport equation as $f$, namely
$$
(\partial_t +  v\cdot \nabla_x + K[f] \cdot \nabla_v) h = 0.
$$
Since the Vlasov flow is divergence-free, a standard energy estimate shows that under the initial assumption $\|\nabla_{x,v} (f_0^\epsilon)\|_{L^\infty(\mathbb{D})}<\infty$, the quantity $\|\nabla_{x,v} ( f(t)^\epsilon)\|_{L^\infty(\mathbb{D})}$ remains finite as long as the Vlasov solution $f(t)$ exists. Therefore, the quantity
\begin{equation}\label{defmaxPT}
P_{T,\epsilon} = \sup_{t\in [0,T]} \left(  \sup\Big\{ |x|+|v|~:~ f(t,x,v) \ne 0\Big\} + \|f(t)\|_{L^\infty} + \|\nabla_v ( f(t)^{\epsilon})\|_{L^\infty(\mathbb{D})} \right)
\end{equation}
is finite for some $T>0$ (we may take $T$ arbitrarily large in dimensions $d\in \{2,3\}$).

Our main result is the validity of the nonlinear Vlasov equation \eqref{VP} as a rigorous limit of the Liouville equation \eqref{NVlasov}. 

\begin{theorem}[Marginal validity of Vlasov equation]\label{theo-main} Let $d\ge 2$ and let $K$ satisfy \eqref{Riesz} with a constant $0\le \alpha <(d-1)/2$. Let $f_0 \ge 0$ be a sufficiently smooth, compactly supported function in $\RR^d\times \RR^d$ such that $\iint f_0 \; dxdv=1$ and $\|\nabla_{x,v} (f_0^{\epsilon})\|_{L^\infty(\mathbb{D})}<\infty$ for some parameter $\epsilon \in (0,1)$ sufficiently small (depending only  on $\alpha$). Let $T>0$ be any fixed time such that the nonlinear Vlasov  equation \eqref{VP} with initial data $f(0,x,v) = f_0(x,v)$ admits a  smooth, compactly supported solution $f(t,x,v)$ on $t \in [0,T]$. Let $F_N(t,z_{[N]})$ be the $N$-particle density distribution function, solving  the Liouville equation  \eqref{NVlasov}, with factorized initial data $F_N(0,z_{[N]}) = f_0^{\otimes N}$, and with $\nu_N = N^{-\delta}$ for some constant $\delta\in (1/2,1)$. Then, for any $k\ge 1$, the $k^{th}$-marginals $\{F_{N,k}\}_{N\ge k}$ satisfy 
\begin{equation}\label{mainconv}
F_{N,k} (t, z_{[k]})\rightharpoonup f(t,z_{[k]})^{\otimes k}
\end{equation}
weakly in $L^2([0,T]\times \DD^k)$ as $N\to \infty$.  
\end{theorem}

Theorem \ref{theo-main} establishes the propagation of the molecular chaos assumption as long as the Vlasov equation is well-posed in all dimensions $d\ge 2$ (which holds up to an arbitrarily long time in dimensions $d\in \{2,3\}$). 

Focusing on the three-dimensional case, our result extends previous works \cite{Jabin,Jabin2} that treat $L^2_{\mathrm{loc}}$ interaction kernels $K(x)$ (thus requiring $\alpha <1/2$). This is the first such result that allows the full range of sub-Coulomb interactions in three dimensions.
Treating the Coulomb interaction, and thus deriving the classical Vlasov-Poisson system in three dimensions, remains an outstanding open problem. While we do not assume any specific sign for the interaction potential in the present paper, the repulsive property may be important for the Coulomb case. In two dimensions, our result covers the Coulomb interaction as well as a class of more singular potentials, which is of independent interest. 
We refer to \cite{Jabin3} for a recent work covering the two-dimensional Coulomb case, both with and without noise.

Due to technical issues, our analysis requires the regularizing effect of the viscosity term coupled with $\nu_N >0$ in \eqref{Nsys}, which helps to absorb the accumulation of lower-order fluctuations. In the main theorem, we choose a small noise with $\nu_N = N^{-\delta}$ with $0<\delta<1$ for simplicity, but the method is flexible enough to cover the case $0<\delta<1$, and also the case $\nu_N\sim 1$ which corresponds to the Vlasov--Fokker--Planck dynamics in the limit. Removing the  microscopic noise to treat the purely deterministic system ($\nu_N = 0$) in three dimensions remains an important open question. 

Finally, note that Theorem \ref{theo-main} establishes the mean-field limit without any explicit rate of convergence. This is due to the fact that our energy-based framework yields uniform bounds on cumulants, from which the BBGKY hierarchy converges weakly to a limiting hierarchy, and the main result follows from the uniqueness of solutions to the limiting hierarchy. Getting quantitative estimates is also an interesting open problem.

\subsection{Outline of the proof}

The proof of Theorem \ref{theo-main} is built upon a weighted energy framework formulated within the dual BBGKY hierarchy. As the first step, following the duality approach introduced by Bresch, Duerinckx, and Jabin \cite{Jabin}, we translate the forward many-body problem into a backward dual equation. More precisely, to estimate the weak convergence of the marginals, we consider a symmetric test function $G_N \sim \psi^{\otimes k}$ and study the dual observable $\Phi_N$ solving the backward equation (see Proposition \ref{prop:dual-problem})
\begin{equation}\label{dualPhiN}
\left\{\begin{aligned}
(\partial_t + \cL^*_N) \Phi_N &= G_N, \qquad \mbox{in}\quad (0,T)\times \DD^N, 
\\
\Phi_N(T) &=0,
\end{aligned}\right.
\end{equation}
where $\cL_N^*$ is the adjoint of the Liouville operator $\cL_N$ given in \eqref{def-LN}. 
Up to a change of variables $(t,x_{[N]},v_{[N]})\mapsto (T-t, x_{[N]}, -v_{[N]})$, we have 
\begin{equation}\label{dualN}
\left\{\begin{aligned}
(\partial_t + \cL_N) \Phi_N &= G_N,
\\
\Phi_N(0) &=0,
\end{aligned}\right.
\end{equation}
for sources of the form  
$$G_N = {{N}\choose{k}}^{-1}\sum_{\sigma_k} \psi^{\otimes k}(t,z_{\sigma_k}).$$

\bigskip
\noindent {\bf Hierarchy of cumulants.} 
As proposed in  \cite{Jabin}, the above duality allows us to express the factorization error exactly in terms of the macroscopic profile $f$ and the dual observable $\Phi_N$ as
\begin{equation}\label{convFNfk}
\langle F_{N,k} - f^{\otimes k}, \psi^{\otimes k}\rangle \approx N\langle (V_{1,2} - \nu_N B_1)\Phi_N, f^{\otimes N}\rangle,
\end{equation}
where 
$$V_{1,2} = (K(x_1-x_2) - K[f](x_1))\cdot \nabla_{v_1} \log f(z_1)$$
captures the force fluctuation and $B_1 = \frac{1}{f(z_1)} \Delta_{v_1}f(z_1)$. To rigorously show that this error vanishes as $N \to \infty$, we decompose $\Phi_N$ into a hierarchy of normalized cumulants  
defined as 
\begin{align}\label{defCNn}
C_{N,n} = \binom{N}{n}^{1/2} \prod_{j=1}^n (1- \mathbb{P}_j) \Phi_{N,n}.
\end{align}
where 
\begin{equation}\label{defPhiNn} \Phi_{N,n} = \int_{\DD^{N-n}} \Phi_N  \; d\mu_{n+1}...d\mu_N= \int_{\DD^{N-n}} \Phi_N  f^{\otimes N-n} \; dz_{n+1}...dz_N
\end{equation}
is the $n$-particle marginal of $\Phi_N$ and $\mathbb{P}_j g_n = \int g_n d\mu_j$ is  the conditional expectation with respect to $d\mu_j = f(t, z_j) dz_j$  for $1\le j\le n$. Then, following \cite[Lemma 11]{Jabin}, the cumulants $C_n=C_{N,n}$ satisfy the following equations
\begin{equation}\label{eqsCNn}
(\partial_t + L_{n} )C_{n} = g_{n}  + T_{n,-}C_{n-1} + T_{n,0}C_{n} + T_{n,1}C_{n+1} + T_{n,2}C_{n+2} 
\end{equation}
with zero initial data $C_{n} (0)=0$ and 
$$g_{n} = {{N}\choose{n}}^{1/2} \prod_{j=1}^n (1- \PP_j) G_{N,n},$$
where $G_{N,n}$ is the n-particle marginal of $G_N$, and 
$$\begin{aligned}
L_n C_n &= \sum_{i=1}^n (v_i \cdot \nabla_{x_i}-\nu_N\Delta_{v_i})C_n + \frac{1}{N-1}\sum_{i\not =j}^n K_{i,j}\cdot \nabla_{v_i} C_n + \frac{N-n}{N-1}\sum_{i=1}^n K[f](x_i)\cdot \nabla_{v_i} C_n,
\\
T_{n,0} C_{n} & = \frac{N-n}{N-1} \sum_{j=1}^n \int_{\DD} V_{*,j}C_n(z_{[n]\setminus\{j\}}, z_*)\; d\mu_* 
+ \frac{1}{N-1}\sum_{i\not =j}^n \nabla_{v_i}\cdot \int_{\DD} K_{i,*} C_n(z_{[n]\setminus\{j\}}, z_*)\; d\mu_*
\\&\quad - \frac{1}{N-1}\sum_{i\not =j}^n  \int_{\DD} V_{*,j} C_n(z_{[n]\setminus\{i\}}, z_*)\; d\mu_*
 + \frac{1}{N-1}\sum_{i\not =j}^n  \int_{\DD^2} K_{*,*'} C_n(z_{[n]\setminus\{i,j\}}, z_*, z_*')\; d\mu_*d\mu_*'\end{aligned}
$$ 
$$
\begin{aligned}
T_{n,-}C_{n-1} & = a_{n,-}\sum_{i\not = j}^n (K[f](x_i) - K_{i,j})\cdot \nabla_{v_i}C_{n-1}(z_{[n]\setminus\{j\}}) 
- a_{n,-} \sum_{i\not =j}^n \int_{\DD} V_{*,j} C_{n-1}(z_{[n]\setminus\{i,j\}}, z_*)\; d\mu_{*}
\\
T_{n,1} C_{n+1} & = 
a_{n,1}
\sum_{j=1}^n 
\int_{\DD} V_{*,j}C_{n+1}(z_{[n]}, z_*, )
\; d\mu_* - a_{n,1}
\sum_{i=1}^n 
\nabla_{v_i}\cdot \int_{\DD} K_{i,*}C_{n+1}(z_{[n]}, z_*, )
\; d\mu_*
\\&\quad -2a_{n,1} \sum_{i=1}^n\int_{\DD^2}V_{*,*'} C_{n+1}(z_{[n]\setminus\{i\}}, z_*, z_*') d\mu_* d\mu_*' + \nu_{n,1} \int_{\DD} B_{n+1}C_{n+1}(z_{[n]}, z_*, )
\; d\mu_*
\\
 T_{n,2}C_{n+2} &=
a_{n,2}
\int_{\DD^2} V_{*,*'} C_{n+2}(z_{[n]}, z_*, z_*')
\; d\mu_*d\mu_{*}'.
\end{aligned}
$$
Here, we have set $a_{n,-}= \frac{\sqrt{N-n+1}}{\sqrt{n}(N-1)}$, $a_{n,1}= \frac{\sqrt{(N-n)(n+1)}}{N-1}$, $a_{n,2} = \frac{\sqrt{(N-n)(N-n-1)(n+1)(n+2)}}{N-1}$, and $\nu_{n,1} = \nu_N \sqrt{(N-n)(n+1)}$.

In the limit of $N\to \infty$ and $\nu_N\to 0$,  we formally have
\begin{equation}\label{expopLn}L_n C_n = \cL^{\otimes n} C_n+ \frac{1}{N-1}\sum_{i\not =j}^n (K_{i,j} - K[f](x_i))\cdot \nabla_{v_i} C_n - \nu_N \sum_{i=1}^n \Delta_{v_i}C_n \approx \cL^{\otimes n} C_n
\end{equation}
where $\cL$  is the one-particle Vlasov operator as in \eqref{VP}. In addition, we note that both $a_{n,-}$ and $a_{n,1}$ tend to zero as $N\to \infty$, and $a_{n,2}\le 2n$. Therefore, it is natural to expect that if 
$$
C_n=C_{N,n} \rightharpoonup  C_{\infty,n}
$$
in a suitable sense, then the limiting cumulants $C_{\infty,n}$ solve
\begin{equation}\label{limitCn}
(\partial_t + \cL^{\otimes n}) C_{\infty,n} = T_{n,0}C_{\infty,n} + T_{n,2}C_{\infty,n+2}
\end{equation}   
with zero initial data $C_{\infty,n}(0)=0$, where
\begin{equation}\label{limitTn}
\begin{aligned}T_{n,0}C_{\infty,n} &= \sum_{j=1}^n \int_{\DD} V_{*,j}C_{\infty,n}(z_{[n]\setminus\{j\}}, z_*)\; d\mu_* ,  
\\T_{n,2}C_{\infty,n+2} &=\sqrt{(n+1)(n+2)}\iint_{\DD^2} V_{*,*'} C_{\infty,n+2}(z_{[n]}, z_*, z_*')
\; d\mu_*d\mu_{*}'.
\end{aligned}
\end{equation}

In this framework, the mean-field limit in $\eqref{mainconv}$ is justified in {\em two main steps:} (1) establishing the convergence $C_{N,n}\rightharpoonup C_{\infty,n}$ as $N\to\infty$ in a suitable sense, and showing that the limiting cumulants  satisfy the hierarchy $\eqref{limitCn}$; and (2) proving uniqueness of solutions to $\eqref{limitCn}$, namely $C_{\infty,n}=0$ for all $t\in[0,T]$. This would justify the convergence of \eqref{mainconv}, upon recalling 
\eqref{convFNfk} and noting $N\langle V_{1,2}\Phi_N, f^{\otimes N}\rangle \approx \langle V_{1,2}, C_{N,2}\rangle$, which vanishes in the limit. 

Concerning the first step, a key observation in  \cite{Jabin} is that thanks to  the orthogonal decomposition \eqref{defCNn}, the cumulants satisfy the following uniform $L^2$-bound
\begin{align}\label{eq:CNn}
\sum_{n\le N} \|C_{N,n}\|_{L^2(\DD^n)}^2 \le C_0.
\end{align}
Consequently, up to a subsequence, for every $n\ge 1$, the cumulants $C_{N,n}$ converges to a weak $L^2$-limit $C_{\infty,n}$ as $N\to \infty$. In case $K \in L^2_{\mathrm{loc}}(\mathbb{D})$, this $L^2$-weak convergence is sufficient to pass from the full dynamics \eqref{eqsCNn} to the limiting equation \eqref{limitCn}. Concerning the second step, as explained in \cite[Lemma 10]{Jabin},  using again the $L^2$-integrability of $K$ and the Cauchy--Schwarz inequality, one may deduce from  \eqref{limitCn} and the zero initial data that 
\begin{align}\label{eq:Cinfty-n-L2-bound}
\|C_{\infty,n}(t)\|_{L^2(\DD^n)}^2 \le C_T n \int_0^t  \Big( \|C_{\infty,n}(s)\|_{L^2(\DD^n)}^2  +  \|C_{\infty,n+2}(s)\|_{L^2(\DD^n)}^2\Big) ds.  
\end{align}
Then the conclusion $C_{\infty,n}=0$ for all time $t\le T$ follows from a Gr\"onwall-type argument, using \eqref{eq:Cinfty-n-L2-bound} and the a priori bound $\sum_{n\ge 0} \|C_{\infty,n}\|_{L^2}^2 \le C_0$, which is a consequence of \eqref{eq:CNn}.  

\medskip

When $K$ is {\em not locally square-integrable}, both steps become significantly harder. Some important progress in this direction has been made recently in \cite{Jabin3} where improved a priori estimates on dual correlations have been obtained by applying hypoelliptic regularity to the hierarchy of equations, a technique inspired by \cite{Bouchut}. 

Our approach below is different from that of \cite{Jabin3}, where we essentially rely on the energy method. To be precise, in the first step we will derive a priori estimates by testing the cumulant equations against a suitable Morawetz-type multiplier; and in the second step we use the Dyson series expansion of the limiting hierarchy and apply the new a priori estimates to conclude. In both steps, the condition $\alpha<(d-1)/2$ in \eqref{Riesz} plays a crucial role. Let us now discuss some details.

\bigskip
\noindent {\bf Morawetz-type estimates.}  
In the present paper, we shall go beyond this $L^2$-barrier on the interaction kernel by devising Morawetz-type weighted estimates to control density concentration at the singularity. Specifically, for the kinetic transport equation $(\partial_t + v\cdot \nabla_x )g=0$, this technique is classical, dating back to the pioneering works of Lions--Perthame \cite{Lions} and Perthame \cite{Perthame} (see also \cite{Perthame1}), to use Morawetz multipliers of the type $\varphi= \frac{x\cdot v}{|x|^\theta}$, with $0\le \theta\le1$, to gain a bulk control of $\iint\frac{|v|^2}{|x|^\theta}g \; dxdv$ in terms of $v$-moments of $g$ (the bound degenerates at $\theta=1$). This provides bounds on density concentration at the origin. We adapt this idea for many-body particle systems. Precisely, for the $n$-particle transport equations of cumulants, we shall use the following Morawetz-type multipliers:
$$
\varphi_n = \sum_{k\not =\ell}^n \frac{(x_k-x_\ell)\cdot (v_k-v_\ell)}{|x_k-x_\ell|^\theta}
$$
for an appropriate choice of $\theta \in (0,1)$. By exploiting the dispersive nature of the free transport operator, the commutator between the transport flow and our Morawetz multiplier generates a positive-definite term that is highly localized near the spatial singularity, namely 
$$\Big(\partial_t + \sum_{i=1}^n v_i \cdot \nabla_{x_i}\Big) \varphi_n \ge (1-\theta) w_n - n^2,\quad w_n = \sum_{k\not =\ell}^n \left( 1+ \frac{|v_k-v_\ell|^2}{|x_k-x_\ell|^\theta} \right).
$$

The remaining force terms in $L_n^*$ are treated as perturbations   
in the resulting weighted energy estimate. This allows to derive an upper bound on $\| \sqrt{w_n}C_n\|_{L^2}$, and the singular weight $w_n$ can be used to control various terms under a weaker condition on $K$. For instance, with $T_{n,2}$ defined as in \eqref{limitCn}, by the Cauchy--Schwarz inequality we have the pointwise estimate 
$$ \begin{aligned}
|T_{n,2}C_{n+2} |
& \le 
a_{n,2}
\int_{\DD^2} |V_{*,*'} C_{n+2}|
\; d\mu_*d\mu_{*}'
\\
& \le 
2n \Big(
\int_{\DD^2}  \left( 1+\frac{|v_{*,*'}|^2}{|x_{*,*'}|^\theta} \right)^{-1} |V_{*,*'}|^2\; d\mu_*d\mu_{*}'\Big)^{1/2} \Big(  \int_{\DD^2}
\left( 1+\frac{|v_{*,*'}|^2}{|x_{*,*'}|^\theta} \right) |C_{n+2}|^2
\; d\mu_*d\mu_{*}'\Big )^{1/2}.
\end{aligned}
 $$
Recall that $V_{*,*'} = (K_{*,*'} - K[f](x_*))\cdot \nabla_{v_*} \log f(z_*)$ where the most singular part is $|K_{*,*'}| \lesssim |x_{*,*'}|^{-\alpha -1}$. Under the assumption $\alpha< (d-1)/2$, we can choose $\max\{0, 2\alpha+2-d\}<\theta<1$ so that 
\begin{equation}\label{intro-bdVstar}\int_{\DD^2} \left( 1+\frac{|v_{*,*'}|^2}{|x_{*,*'}|^\theta} \right)^{-1} |V_{*,*'}|^2\; d\mu_*d\mu_{*}' \le C_T. 
\end{equation}
Here to prove \eqref{intro-bdVstar} we use 
$$\left( 1+\frac{|v_{*,*'}|^2}{|x_{*,*'}|^\theta} \right)^{-1} |K_{*,*'}|^2 \lesssim \left( \frac{|x_{*,*'}|^\theta}{|v_{*,*'}|^2}\right)^{1-\epsilon}|x_{*,*'}|^{-2\alpha-2} =  \frac{ |x_{*,*'}|^{\theta(1-\epsilon)-2\alpha-2} }{ |v_{*,*'}|^{2(1-\epsilon)} } $$
which is  locally integrable. On the technical side,  $V_{*,*'}$ is slightly more singular than $K_{*,*'}$ since $\nabla_{v} \log f$ is not uniformly bounded, but we can use $f^\epsilon |\nabla_v \log f|\lesssim_{\epsilon} |\nabla_v (f^\epsilon)| \lesssim_{\epsilon,T} 1$ from the assumption \eqref{defmaxPT} to deduce that 
\begin{equation}\label{eq:V12-upper-bound}
|V_{*,*'}|^2 d\mu_* d\mu_{*'} \lesssim_{\epsilon,T} ( |x_{*,*'}|^{-2\alpha-2}+1)f_{*}^{1-2\epsilon }f_{*'} dz_{*} dz_{*'}, 
\end{equation}
which is sufficient to ensure \eqref{intro-bdVstar}. Combining with the symmetry of $C_{n+2}$, we obtain
\begin{align}\label{eq:Tn+2-intro}
\int_{\DD^n}|T_{n,2}C_{n+2} |^2 d \mu_1 ... d\mu_n \le C_T  \int_{\DD^{n+2}} w_{n+2} 
 |C_{n+2}|^2
\; d\mu_1 ... d\mu_{n+2}.
\end{align}

In this way, we can close the estimates using the weighted norm $\|\sqrt{w_n} C_n\|_{L^2_{t,z}}$. On the technical side, the absorption of lower-order errors also requires the regularizing properties of the microscopic viscosity, which is the reason we need to add some noise to the model. To be precise, after an iteration procedure, we obtain the a priori estimate  
\begin{equation}\label{key-intro}
\int_0^T \| \sqrt{w_n} C_n\|_{L^2_{z}}^2 dt  \le C_{T} n^{b}
\end{equation}
for a large but fixed constant $b$. Using \eqref{key-intro} we can leverage the weak compactness framework from \cite{Jabin} to our weighted $L^2$-setting, and  obtain the limiting hierarchy \eqref{limitCn}.

\bigskip
\noindent {\bf Uniqueness of limiting hierarchy.} Another key ingredient of our proof is the uniqueness of the solution of \eqref{limitCn}. The main difficulty is that \eqref{eq:Cinfty-n-L2-bound} is unavailable since the force kernel $K$ needs not belong to $L^2_{\mathrm{loc}}$. 

In order to use \eqref{key-intro}, we will reformulate \eqref{limitCn} in terms of the Dyson series, using the key observation that time integration along the Vlasov characteristics produces a regularizing effect. Here the symmetry of cumulants is helpful to reorganize the time-ordered Dyson expansion into products of time-integrated two-particle operators. To be precise, we denote by $U^{(j)}$ the propagator generated by 
\begin{equation}\label{cLparticlej-intro}
	\cL^{(j)} C_n = v_j\cdot \nabla_{x_j}C_n+ K[f](x_j)\cdot \nabla_{v_j}C_n  - \int_{\DD} V_{*,j}C_n(z_{[n]\setminus\{j\}}, z_*)\; d\mu_*.
\end{equation}
which contains the one-particle Vlasov
operator  acting on the variable $z_j$, plus 
a bounded nonlocal correction arising from the
limiting hierarchy.  Then by iterating Duhamel's formula from \eqref{limitCn}, for every $k\ge 1$, we can write 
\begin{align*}
	C_{\infty,n}(t) & =  \int_0^t U_n(t,s)  T_{n,2}C_{\infty,n+2}(s)\; ds  \\
	&=  \frac{1}{k!} \sqrt{\frac{(n+2k+2)!}{n!}}  \int_0^t \prod_{j=1}^n U^{(j)}(t,s) 
			\prod_{j=1}^k\Big(\int_s^t\cB_{n+2j}(\tau,s) \; d\tau\Big) T_{n+2k}(s) C_{\infty,n+2k+2}(s)\; ds
\end{align*}
where $T_{n,2}= \sqrt{(n+1)(n+2)} T_{n+2}$ with  $T_{n,2}$ given in \eqref{limitCn} and 
	$$
	\cB_{n+2j}(t,s) = T_{n+2j}(t) U^{(n+2j-1)}(t,s) U^{(n+2j)} (t,s). 
	$$
	All operators $\cB_{n+2j}$ appearing in the Dyson series have the same structure which is captured by the  two-particle operator 
	\begin{equation}\label{eq:B-intro}
\cB(t,s)=T_2(t) U^{(1)}(t,s) U^{(2)}(t,s). 
\end{equation}
Here the main difficulty lies in the singular interaction associated with the operator
$$
T_2(t) G
=\int_{\DD^2}V_{1,2}(t) G(z_1,z_2)\,d\mu_1d\mu_2
$$
where $V_{1,2}(t) = (K_{1,2} - K[f](t,x_1))\cdot \nabla_{v_1} \log f(t,z_1)$. As we mentioned before, if $K \notin L^2_{\mathrm{loc}}$, then $T_2$ is not bounded in  $L^2$, and the same applies to $\cB$ in \eqref{eq:B-intro}. However, our key observation is that in the Dyson series, the time integration gains regularity through phase mixing, namely 
\begin{equation}\label{eq:B-bound-intro}
\left|\int_s^t \cB(\tau,s)G\,d\tau\right|
\lesssim |t-s|^{1-\theta/2}\|G\|_{L^2(\DD^2)}.
\end{equation}
See Lemma \ref{lem-bdCBt}. 

Roughly speaking, the nonlocal correction in $\cL^{(j)}$ is bounded in $L^2$, and hence it can be treated perturbatively by  Duhamel's formula. Thus the regularizing effect essentially comes from the pure Vlasov characteristics. By changing the variables and  the uniform bound of the corresponding Jacobian, we reduce the problem to that of the free transport flow $(x_j,v_j) \mapsto (x_j -tv_j,v_j)$, where the relative displacement of two particles at different times $\tau_1,\tau_2$ is proportional to $y=(\tau_1-\tau_2)(v_1-v_2)$. In this way, we deduce \eqref{eq:B-bound-intro} from a ``$TT^*$ argument" and the convolution-type estimate
\begin{equation}\label{eq:K-K-bound-intro}
\int_{|x|\lesssim 1}|K(x)|\,|K(x-y)|\,dx
\lesssim 1+|y|^{-\theta}
\end{equation}
for any constant $\max\{0,2\alpha+2-d\}<\theta<1$. Here the lower bound $\theta>2\alpha+2-d$ comes from the singularity $|K(x)|\lesssim |x|^{-1-\alpha}$, and the upper bound $\theta<1$ is taken to ensure the time integrability.  Again the condition $\alpha<(d-1)/2$ is crucial to make sure that these requirements are compatible. 

In summary, by inserting the bound \eqref{eq:B-bound-intro} into the time-ordered Dyson expansion, we find that each time iteration produces a small factor proportional to $t^{1-\theta/2}$. Together with the a priori bound \eqref{key-intro} on the cumulants, after $k$ iteration steps, we obtain 
an estimate of the form
\begin{equation}
\|C_{\infty,n}(t)\|_{L^2(\DD^n)}  \le (C_{T,n} k)^{(n+b)/2+1}(C_Tt^{1-\theta/2})^k,
\end{equation}
which vanishes as $k\to\infty$ for $t\le t_*$ if $t_*$ is chosen sufficiently small such that $C_Tt_*^{1-\theta/2} <1$. By a standard continuation argument, we conclude the uniqueness of the limiting hierarchy for all time $t\in [0,T]$ and complete the proof.

\bigskip
\noindent {\bf Organization of the paper.} In Section \ref{sec:Duality approach} we quickly recall the duality framework and the formulation of the cumulant dynamics. In Section \ref{sec-linear}, we establish the key linear weighted estimate based on the Morawetz-type multiplier. Then in Section \ref{sec:Iterative bounds} we explain the iterative bootstrap procedure for the cumulant hierarchy, where the absorption of lower-order errors requires the regularizing properties of the microscopic viscosity. Then in Section \ref{sec:final-cumulant} we derive uniform bounds on the full cumulant hierarchy. Finally, in Section \ref{sec:final-uniqueness}, we pass $C_{N,n} \rightharpoonup  C_{\infty,n}$ as $N\to \infty$ and justify the uniqueness of the limiting hierarchy, thus concluding the proof of Theorem \ref{theo-main}.

\bigskip
\noindent{\bf Notations.} Throughout the paper, we often drop the dependence on $N$ in the cumulants $C_{N,n}$ and simply write $C_n$ when there is no ambiguity. The notation $A \lesssim B$ means there exists a constant $c>0$ independent of $N,n$ such that $A \le c B$.

\bigskip
\noindent
{\bf Acknowledgments.} We would like to thank Pierre-Emmanuel Jabin and Mitia Duerinckx for inspiring discussions. We thank the Vietnam Institute for Advanced Study in Mathematics and its staff for their warm hospitality. TN also acknowledges the hospitality of the Institut des Hautes \'Etudes Scientifiques during his visits. Part of this work was carried out at both institutes. PTN acknowledges the partial support from the European Research Council via the ERC Consolidator Grant RAMBAS, Project Nr. 10104424. TN is supported in part by the NSF under grant DMS-2349981 and a Simons Foundation fellowship. 

\section{Duality approach} \label{sec:Duality approach}

In this section, we quickly review the duality approach introduced by Bresch, Duerinckx, and Jabin in \cite{Jabin}. For the reader's convenience, we provide some elementary proofs and fix our notation. 

\subsection{Dual problem}  The starting point is the following simple proposition which is a direct consequence of \eqref{NVlasov}.

\begin{proposition} \label{prop:dual-problem} Let $G_N\in C_0^\infty(\RR\times \DD^N)$ be any symmetric smooth test function, and let 
$\Phi_N$ be any solution to the dual problem 
\begin{equation}\label{dualPhiN}
\left\{\begin{aligned}
(\partial_t + \cL^*_N) \Phi_N &= G_N, \qquad \mbox{in}\quad (0,T)\times \DD^N, 
\\
\Phi_N(T) &=0,
\end{aligned}\right.
\end{equation}
for some fixed $T>0$, where $\cL_N^*$ is dual to $\cL_N$ with respect to the usual $L^2$ inner product in $\DD^N$, namely 
$$\cL_N^* = \sum_{i=1}^N (v_i \cdot \nabla_{x_i} +\nu_N \Delta_{v_i})+  \frac{1}{N-1}\sum^N_{i\not =j} K(x_i - x_j)\cdot \nabla_{v_i}
.$$ 
Then, there holds 
\begin{equation}\label{dual-form}
\begin{aligned}
\langle F_N - f^{\otimes N}, G_N \rangle_{L^2([0,T]\times \DD^N)}  &= N\langle (V_{1,2} - \nu_N B_1)\Phi_N, f^{\otimes N}\rangle_{L^2([0,T]\times \DD^N)}
\\&\quad - \langle F_N(0) - f_0^{\otimes N}, \Phi_N(0) \rangle_{L^2(\DD^N)},
\end{aligned}
\end{equation}
in which $V_{1,2} = (K(x_1-x_2) - K[f](x_1))\cdot \nabla_{v_1} \log f(z_1)$ and $B_1 = \frac{1}{f(z_1)} \Delta_{v_1}f(z_1)$.
\end{proposition}

\begin{proof} For any symmetric test functions $\Phi_N(t,z_1,\cdots, z_N)$, using  \eqref{NVlasov}, we compute 
$$\begin{aligned}
\frac{d}{dt}\langle F_N, \Phi_N\rangle 
&= \langle \partial_t F_N, \Phi_N\rangle + \langle F_N, \partial_t \Phi_N\rangle 
\\
&= -\langle \cL_N F_N, \Phi_N\rangle + \langle F_N, \partial_t \Phi_N\rangle 
= \langle F_N, (\partial_t +\cL^*_N)\Phi_N\rangle,
\end{aligned}
$$
upon noting that $\cL_N$ can be written in the divergence form. Here and in what follows, 
$\langle f,g \rangle = \int f(z)g(z)\; dz$, the usual $L^2$ inner product in $\DD^N$.  
Similarly, we compute 
$$\begin{aligned}
\frac{d}{dt}\langle f^{\otimes N}, \Phi_N\rangle 
&= \langle f^{\otimes N}, (\partial_t +\cL^{\otimes N})\Phi_N\rangle 
\end{aligned}
$$
where $\cL$ is the one-particle operator as defined in \eqref{VP}, and 
$$ \cL^{\otimes N} = \sum_{i=1}^N v_i \cdot \nabla_{x_i}+ \sum_{i=1}^N K[f](x_i)\cdot \nabla_{v_i}.$$
Therefore, we obtain 
\begin{equation}\label{dual-id}
\begin{aligned}
\frac{d}{dt}\langle F_N-f^{\otimes N}, \Phi_N\rangle 
&= \langle F_N, (\partial_t +\cL_N^*)\Phi_N\rangle - \langle f^{\otimes N}, (\partial_t +\cL^{\otimes N})\Phi_N\rangle 
\\
&= \langle F_N - f^{\otimes N}, (\partial_t +\cL_N^*)\Phi_N\rangle + \langle f^{\otimes N}, (\cL_N^* - \cL^{\otimes N})\Phi_N\rangle 
.\end{aligned}\end{equation}
By definition, we note 
$$
\begin{aligned}
&\langle f^{\otimes N}, (\cL_N^* - \cL^{\otimes N})\Phi_N\rangle 
\\&=\frac{1}{N-1} \sum_{i\not =j} \langle f^{\otimes N},\Big[ K(x_i-x_j) - K[f](x_i)\Big] \cdot \nabla_{v_i} \Phi_N\rangle + \nu_N \sum_{i=1}^N\langle f^{\otimes N}, \Delta_{v_i}\Phi_N\rangle 
\\& = - \frac{1}{N-1} \sum^N_{i\not =j} \langle f^{\otimes N},V_{i,j} \Phi_N\rangle + \nu_N \sum_{i=1}^N\langle f^{\otimes N}, B_i\Phi_N\rangle,
\end{aligned}$$
in which $V_{i,j} = (K(x_i-x_j) - K[f](x_i))\cdot \nabla_{v_i} \log f(z_i)$ and $B_i = \frac{1}{f(z_i)}\Delta_{v_i}f(z_i)$. In addition, by the symmetry of $\Phi_N$,  we can further write 
$$
\begin{aligned}
\langle f^{\otimes N}, (\cL_N^* - \cL^{\otimes N})\Phi_N\rangle 
& = - N\langle (V_{1,2} - \nu_N B_1) \Phi_N, f^{\otimes N}\rangle .
\end{aligned}$$
This leads to 
\begin{equation}\label{dual-id1}
\begin{aligned}
\frac{d}{dt}\langle F_N-f^{\otimes N}, \Phi_N\rangle 
&= \langle F_N - f^{\otimes N}, (\partial_t +\cL_N^*)\Phi_N\rangle  - N\langle (V_{1,2} - \nu_N B_1)\Phi_N, f^{\otimes N}\rangle
\end{aligned}\end{equation}
for any symmetric test function $\Phi_N$. Integrating this in time and taking $\Phi_N$ solving the dual problem \eqref{dualPhiN}, we obtain the proposition. 
\end{proof}

\begin{corollary}\label{cor-main}
Fix $k\ge 1$, and take 
\begin{equation}\label{sourceGN}
G_N(t,z_{[N]}) = {{N}\choose{k}}^{-1}\sum_{\sigma_k} \psi^{\otimes k}(t,z_{\sigma_k})
\end{equation}
for any $\psi \in C_0^\infty([0,T] \times \RR^d\times \RR^d)$, where $\sigma_k = (i_1,\cdots, i_k)$ is any $k$-element subset of $\{1,\cdots,N\}$. Then, 
\begin{equation}\label{dual-formk}
\begin{aligned}
\langle F_{N,k} - f^{\otimes k}, \psi^{\otimes k}\rangle_{L^2([0,T]\times \DD^k)}  &= N\langle (V_{1,2} - \nu_N B_1) \Phi_N, f^{\otimes N}\rangle_{L^2([0,T]\times \DD^N)}
\\&\quad - \langle F_N(0) - f_0^{\otimes N}, \Phi_N(0) \rangle_{L^2(\DD^N)},
\end{aligned}
\end{equation}
for any solution $\Phi_N$ to the dual problem \eqref{dualPhiN}. 
\end{corollary}
\begin{proof} The corollary is direct from the previous proposition, upon noting that 
$$\langle F_N - f^{\otimes N}, G_N \rangle_{L^2([0,T]\times \DD^N)}   = \langle F_{N,k} - f^{\otimes k}, \psi^{\otimes k}\rangle_{L^2([0,T]\times \DD^k)} $$
that follows from the choice of $G_N$. 
\end{proof}

Up to a change of variables $(t,x_{[N]},v_{[N]})\mapsto (T-t, x_{[N]}, -v_{[N]})$, we obtain  \eqref{dualN} 
with the same operator $\cL_N$ defined as in \eqref{def-LN}. We are interested in estimating the right-hand side of \eqref{dual-formk} in the mean-field limit of $N\to \infty$. To proceed, in what follows, we shall use the probability measure $d\mu = fdz$, namely 
$$ \langle F,G\rangle_n = \int_{\DD^n} FG d\mu_1...d\mu_n= \int_{\DD^n} FG f^{\otimes n} \;d{z_{[n]}}.$$
Note that the Vlasov equation \eqref{VP} is invariant under the change of variables $(t,x,v)\mapsto (T-t, x,-v)$, and therefore, $f(t,x,v)$ remains to solve the same Vlasov equation \eqref{VP}. In addition, thanks to the conservation of mass $\iint f(t)\;dxdv=1$, we have $\langle 1\rangle_n =1$ for any $n\ge 1$.

\subsection{Marginal equations}

In this section, we consider the dual problem  \eqref{dualN}.  We use notion of marginals in \eqref{defPhiNn}, and also introduce
the $n$-marginal operator 
\begin{equation}\label{defLNn}
L_{N,n} = \sum_{i=1}^n (v_i \cdot \nabla_{x_i} - \nu_N \Delta_{v_i})+ \frac{1}{N-1}\sum_{i\not =j}^n K(x_i-x_j) \cdot \nabla_{v_i} 
.
\end{equation}
Note that the marginal operator $L_{N,n}$ in \eqref{defLNn} differs from the cumulant operator $L_n$ in \eqref{expopLn} since it does not include the self-consistent mean-field force (we will come to $L_n$ in the next subsection). 

\begin{proposition}\label{prop-PhiNn} Let $\Phi_N$ be any solution to \eqref{dualN}, and $\Phi_{N,n}$ be its $n$-particle marginals defined as in \eqref{defPhiNn}. 
Then, for any $1\le n\le N-2$, we have
\begin{equation}\label{eqsPhiNn}
(\partial_t + L_{N,n})\Phi_{N,n} = G_{N,n} + D_{N,n}^{(1)}\Phi_{N,n+1} + S_{N,n}^{(1)} \Phi_{N,n+1} + S_{N,n}^{(2)}\Phi_{N,n+2} 
\end{equation}
with zero initial data $\Phi_{N,n} (0)=0$, where 
$$
\begin{aligned}
D_{N,n}^{(1)}\Phi_{N,n+1} & =  -\frac{N-n}{N-1}\sum_{i=1}^n\nabla_{v_i} \cdot
\int_{\DD} K_{i,n+1}
\Phi_{N, n+1}\; d\mu_{n+1}
\\
S_{N,n}^{(1)} \Phi_{N,n+1} & = 
\frac{N-n}{N-1}\sum_{i=1}^n 
\int_{\DD} \Phi_{N,n+1}V_{n+1,i}
\; d\mu_{n+1} +\nu_N (N-n)\int_{\DD} B_{n+1}\Phi_{N,n+1} \; d\mu_{n+1}
\\
 S_{N,n}^{(2)}\Phi_{N,n+2} &=
\frac{(N-n)(N-n-1)}{N-1}
\int_{\DD^2} \Phi_{N,n+2}V_{n+1,n+2} 
\; d\mu_{n+1}d\mu_{n+2}.
\end{aligned}
$$
\end{proposition}

\begin{proof} Introduce $L_{N,n}$ and $\Phi_{N,n}$ as in \eqref{defPhiNn}-\eqref{defLNn}. Integrating \eqref{dualN} over $\DD^{N-n}$ with volume elements $d\mu_{n+1}...d\mu_N$, we obtain 
$$
\begin{aligned} 
G_{N,n} 
&= \int_{\DD^{N-n}} (\partial_t + \cL_N)\Phi_N
\\
&= \int_{\DD^{N-n}} (\partial_t + L_{N,n})\Phi_N + \int_{\DD^{N-n}} (\cL_N -  L_{N,n})\Phi_N.
\end{aligned}$$ 
Using $\partial_t f^{\otimes (N-n)} /f^{\otimes (N-n)} = \sum_{i\ge n+1} \log f(z_i)$, we compute 
$$
\begin{aligned}
 \int_{\DD^{N-n}} (\partial_t + L_{N,n})\Phi_N 
 &= (\partial_t + L_{N,n})\Phi_{N,n} - \sum_{i\ge n+1} \int_{\DD^{N-n}} \Phi_N \partial_t\log f(z_i).
 \end{aligned}$$
 On the other hand, by the definition of $L_{N,n}$, we have 
$$ \cL_N -  L_{N,n} = \sum_{i\ge n+1} (v_i \cdot \nabla_{x_i} -\nu_N  \Delta_{v_i})+ \frac{1}{N-1}\Big( \sum_{i\le n, \; j\ge n+1} + \sum_{j\le n, \; i\ge n+1} + \sum_{i\not =j\ge n+1}\Big ) K_{i,j}\cdot \nabla_{v_i}$$
in which $K_{i,j} = K(x_i-x_j)$.
Therefore, upon integrating by parts in $z_i$, we obtain  
\begin{equation}\label{compPNn}
\begin{aligned}
(\partial_t + L_{N,n})\Phi_{N,n} 
&= G_{N,n} + \sum_{i\ge n+1}\int_{\DD^{N-n}} \Phi_N \Big[(\partial_t+ v_i\cdot \nabla_{x_i}) \log f(z_i) + \nu_N \frac{1}{f(z_i)}\Delta_{v_i}f(z_i)\Big]\\&\quad - \frac{1}{N-1}\Big( \sum_{i\le n, \; j\ge n+1} + \sum_{j\le n, \; i\ge n+1} + \sum_{i\not =j\ge n+1}\Big )\int_{\DD^{N-n}} K_{i,j}\cdot \nabla_{v_i} 
\Phi_{N}.
\end{aligned}
\end{equation}
Let us compute each term on the right. Indeed, recalling that $f$ solves \eqref{VP} and using the fact that $\Phi_N$ is symmetric in $z_{[N]}$, we have 
$$
\begin{aligned} 
&\sum_{i\ge n+1}\int_{\DD^{N-n}} \Phi_N \Big[(\partial_t+ v_i\cdot \nabla_{x_i}) \log f(z_i) + \nu_N \frac{1}{f(z_i)}\Delta_{v_i}f(z_i)\Big]
\\&= - (N-n)\int_{\DD} \Phi_{N,n+1} \Big[ K[f](x_{n+1})\cdot \nabla_{v_{n+1}}\log f(z_{n+1}) - \nu_N B_{n+1}\Big]\; d\mu_{n+1} ,
\end{aligned}
$$
upon recalling $B_i =\frac{1}{f(z_i)}\Delta_{v_i}f(z_i)$. 
Next, we compute 
$$
\begin{aligned}
\sum_{i\le n, \; j\ge n+1}
\int_{\DD^{N-n}} K_{i,j}\cdot \nabla_{v_i} 
\Phi_{N}
& = (N-n)\sum_{i\le n}\nabla_{v_i} \cdot
\int_{\DD} K_{i,n+1}
\Phi_{N, n+1}\; d\mu_{n+1}
\\\sum_{j\le n, \; i\ge n+1}
\int_{\DD^{N-n}} K_{i,j}\cdot \nabla_{v_i} 
\Phi_{N} 
& = -(N-n)\sum_{j\le n} 
\int_{\DD} \Phi_{N,n+1}K_{n+1,j} \cdot\nabla_{v_{n+1}}\log f(z_{n+1})
\; d\mu_{n+1}.
\end{aligned}
$$
Similarly, we have 
$$
\begin{aligned}
&\sum_{i\not =j\ge n+1}
\int_{\DD^{N-n}} K_{i,j}\cdot \nabla_{v_i} 
\Phi_{N}
\\& =-(N-n)(N-n-1)
\int_{\DD^2} \Phi_{N,n+2}K_{n+1,n+2}\cdot \nabla_{v_{n+1}} \log f(z_{n+1}) 
\; d\mu_{n+1}d\mu_{n+2}.
\end{aligned}
$$
Combining, we thus obtain 
\begin{equation}\label{compPNn1}
\begin{aligned}
(\partial_t + L_{N,n})\Phi_{N,n} 
&= G_{N,n} - (N-n)\int_{\DD} \Phi_{N,n+1} \Big[ K[f](x_{n+1})\cdot \nabla_{v_{n+1}}\log f(z_{n+1}) - \nu_N B_{n+1}\Big] \; d\mu_{n+1}
\\&\quad -\frac{N-n}{N-1}\sum_{i\le n}\nabla_{v_i} \cdot
\int_{\DD} K_{i,n+1}
\Phi_{N, n+1}\; d\mu_{n+1}
\\&\quad + \frac{N-n}{N-1}\sum_{i\le n} 
\int_{\DD} \Phi_{N,n+1}K_{n+1,i} \cdot\nabla_{v_{n+1}}\log f(z_{n+1})
\; d\mu_{n+1}
\\&\quad + \frac{(N-n)(N-n-1)}{N-1}
\int_{\DD^2} \Phi_{N,n+2}K_{n+1,n+2}\cdot \nabla_{v_{n+1}} \log f(z_{n+1}) 
\; d\mu_{n+1}d\mu_{n+2}.
\end{aligned}
\end{equation}
Finally, recalling the notation $V_{i,j} = (K_{i,j} - K[f](x_i))\cdot \nabla_{v_i} \log f (z_i)$, we may replace $K_{i,j}\cdot \nabla_{v_i} \log f (z_i)$ by $V_{i,j}$ in the above last two integral terms, creating the following extra terms  
$$
\begin{aligned}
&  \frac{N-n}{N-1}\sum_{i\le n} 
\int_{\DD} \Phi_{N,n+1}K[f](x_{n+1})\cdot\nabla_{v_{n+1}}\log f(z_{n+1})
\; d\mu_{n+1}
\\& + \frac{(N-n)(N-n-1)}{N-1}
\int_{\DD^2} \Phi_{N,n+2}K[f](x_{n+1})\cdot \nabla_{v_{n+1}} \log f(z_{n+1}) 
\; d\mu_{n+1}d\mu_{n+2}
\\&=(N-n)\int_{\DD} \Phi_{N,n+1} K[f](x_{n+1})\cdot \nabla_{v_{n+1}}\log f(z_{n+1}) \; d\mu_{n+1}
\end{aligned}
$$
which cancels out with part of the second term on the right-hand side of \eqref{compPNn1}. This yields the proposition as stated. 
\end{proof}

\subsection{Cumulant equations}

In view of Corollary \ref{cor-main}, our aim is to prove that 
\begin{equation}\label{recallgoal}
\begin{aligned}
&N\langle (V_{1,2} - \nu_N B_1) \Phi_N, f^{\otimes N}\rangle_{L^2([0,T]\times \DD^N)}  
\\&= N \int_0^T \int_{\DD^2} V_{1,2}\Phi_{N,2} \; d\mu_1d\mu_2  - \nu_NN \int_0^T \int_{\DD}B_1\Phi_{N,1} \; d\mu_1 
\to 0
\end{aligned}
\end{equation}
as $N\to \infty$, in which the dynamics of $\Phi_{N,2}$ is driven by $G_{N,2}$ through \eqref{eqsPhiNn} (and of course, being coupled with higher order marginals). Since the source $G_N$ is of the form \eqref{sourceGN}, $G_{N,2}$ may contain nonzero constants. Indeed, for instance, we take $k=1$ and 
$G_N(t,z_{[N]}) = \frac1N\sum_{j=1}^N \psi(t,z_j). 
$
By a direct calculation, we have 
$$G_{N,n} = \frac1N \sum_{j=1}^n \psi(t,z_j) + \frac{N-n}{N}\langle \psi\rangle_\DD,$$
in which $\langle \psi\rangle_\DD$ denotes the average of $\psi$ over the phase space $\DD$. As a result, $G_{N,n}$ does not tend to zero, and neither does $\Phi_{N,n}$, in the limit of $N\to \infty$. Observe however that constants do not contribute into the convergence \eqref{recallgoal}, since $\int B_1d\mu_1=0$ and $\int V_{1,2} d\mu_j =0$ for $j=1$ or $j=2$. We recall the following results from \cite[Lemma 11]{Jabin} on cumulants $C_{N,n}$. 

\begin{proposition}\label{prop-CNn} Let $\Phi_N$ and $\Phi_{N,n}$ be defined as in \eqref{dualN}-\eqref{defPhiNn}. Then the cumulants $C_n=C_{N,n}$ defined in  \eqref{defCNn} satisfy the equation \eqref{eqsCNn}.  
\end{proposition}

In order to study the behavior of $C_{N,n}$ when $N\to \infty$, estimates uniform in $N$ will play a crucial role. The following uniform bound, which is based on the orthogonality of $\{C_{N,n}\}_n$, is a key ingredient in \cite{Jabin}.  

\begin{lemma}\label{lem-apriori}
Let $C_{N,n}$ be the normalized cumulants as in \eqref{defCNn}. Then, for $0\le \ell \le 2$, there hold 
\begin{equation}\label{L2bounds-apriori} \sup_{[0,T]} \sum_n \|C_{N,n}\|_{L^2(\DD^n)}^2 \le C_0, \qquad \sup_{[0,T]} \sum_n \Big(\frac{n}{N}\Big)^\ell \|C_{N,n}\|_{L^\infty(\DD^\ell; L^2(\DD^{n-\ell}))}^2 \le C_0 ,
\end{equation}
for some universal constant $C_0$.
\end{lemma}

New bounds on $C_{N,n}$ will be derived in the next section. In what follows, we shall often drop the apparent dependence on $N$ and write $C_n=C_{N,n}$.

\section{New linear estimates}\label{sec-linear}

In this section we study the linear problem 
\begin{equation}\label{pLn}
(\partial_t + L_{n}) C_{n} = R_{n} 
\end{equation}
on $[0,T]\times \DD^n$, with zero initial data $C_{n}(0) =0$, where $L_n$ is defined as in Proposition \ref{prop-CNn}. We keep the error term $R_n$ as a general source, which will be specified later. Here we consider the $L^2$ spaces with respect to probability measure $d\mu_j = f dz_j$, namely 
\begin{equation}\label{defL2norm} 
\begin{aligned}
\langle C_{n}, D_n \rangle_{L^2(\DD^n)} &= \int_{\DD^n} C_{n} D_n \; d\mu_1\cdots d\mu_n
\end{aligned}
\end{equation}
together with the $L^2$ norm $\|C_n\|_{L^2}= \sqrt{\langle C_n, C_n \rangle_{L^2}}$. We also use the weighted $L^2$ spaces associated with the norms $\|\cdot \|_{L^2_w} = \| \sqrt{w}\cdot \|_{L^2}$ for some weight function $w = w(z_{[n]}) \ge 0$.

In Section \ref{sec:Weighted L2} we derive a weighted $L^2$-estimate on the solution $C_n$ for \eqref{pLn} with a source $R_n\in L^2$, which is the key ingredient of our approach. In Section \ref{sec:H-1sources}, we extend this weighted estimate to a general source $R_n$ in $H^{-1}$, which is helpful for an iteration argument in Section \ref{sec:Iterative bounds}.

\subsection{Weighted $L^2$ estimates}\label{sec:Weighted L2}

The main estimate in this subsection is the following. 

\begin{proposition}\label{prop-keyL2} Let $d\ge 2$. Let $K(x)$ be the Riesz interaction force \eqref{Riesz} with $0\le \alpha < \frac{d-1}{2}$, and let $C_n$ be the solution to \eqref{pLn} and satisfy the a priori bounds in Lemma \ref{lem-apriori}. Then, for $\nu_Nn\lesssim 1$, there are positive constants $\delta_\alpha$ and $C_{\alpha,T}$, independent of $n$, so that 
\begin{equation}
\label{keysupL2}
\begin{aligned}
n^2\sup_{[0,T]}\|C_n(t)\|_{L^{2}(\DD^n)}^2 &+ \nu_Nn^2\sum_i^n\|\nabla_{v_i}C_n\|^2_{L^{2}([0,T]\times \DD^n)}  + \| \sqrt{w_n} C_n\|_{L^{2}([0,T]\times \DD^n)}^2 
\\
&\le 
C_{\alpha,T}n^2 \| R_n\|_{L^{2}([0,T]\times \DD^n)}^2 +
C_{\alpha,T}n^4(N/n)^{-2\delta_\alpha }
\end{aligned}
\end{equation}
for all $1\le n\le N$, where the weight $w_n \ge 0$ is defined by 
\begin{equation}\label{weight}
w_n = \sum_{k\not =\ell}^n \left( 1+ \frac{|v_{k,\ell}|^2}{|x_{k,\ell}|^\theta} \right) 
\end{equation}
for some constant $\theta \in (0,1)$, where $x_{k,\ell} := x_k-x_\ell$ and $v_{k,\ell} := v_k - v_\ell$. 
\end{proposition}

Unlike Lemma \ref{lem-apriori} which provides a priori boundedness of $C_n$, Proposition \ref{prop-keyL2} provides decay and stability estimates in the large $N$ limit. The constraint $\nu_Nn\lesssim1$ is imposed for sake of simplicity in the presentation, and will be removed in the next section. The weighted norms control the spatial singularity at $x_k = x_\ell$, which plays a crucial role in treating the singularity of $K(x)$.

\medskip
Let us start with a basic $L^2$-estimate.

\begin{lemma}\label{lem-L2basic} Let $C_n$ be the solution to \eqref{pLn}. Then, there holds
\begin{equation}\label{L2basic}
\begin{aligned}
&\frac{d}{dt}\|C_n\|_{L^2(\DD^n)}^2 +\nu_N\sum_i^n\|\nabla_{v_i}C_n\|^2_{L^2(\DD^n)}
\\&\le 2  \langle R_n,C_n\rangle_{L^2(\DD^n)}  + \frac{n^2}{N-1} \langle |V_{1,2}|, C_n^2 \rangle_{L^2(\DD^n)} 
 + C_T \sqrt{\nu_N} n,
\end{aligned}\end{equation}
where $C_T$ denotes some constant that depends on $T,\delta$ and the initial data $f_0$, but does not depend on $n, N$.
Here $V_{1,2} = (K(x_1-x_2) - K[f](x_1))\cdot\nabla_{v_1}\log f(z_1)$. 
\end{lemma}
\begin{proof} Recall from Proposition \ref{prop-CNn} that 
\begin{align}\label{opLn}
L_n C_n &= \sum_{i=1}^n (v_i \cdot \nabla_{x_i} + K[f](x_i)\cdot \nabla_{v_i} -\nu_N\Delta_{v_i})C_n + \frac{1}{N-1}\sum_{i\not =j}^n (K_{i,j} - K[f](x_i))\cdot \nabla_{v_i} C_n\nonumber\\
&= \cL^{\otimes n} C_n+ \frac{1}{N-1}\sum_{i\not =j}^n (K_{i,j} - K[f](x_i))\cdot \nabla_{v_i} C_n, 
\end{align}
where $\cL$ denotes the one-particle Vlasov  operator as in \eqref{VP}. Therefore,  integrating the equation \eqref{pLn} against $C_n$ and upon recalling $d\mu_{[n]} = f^{\otimes}dz_{[n]}$, we compute  
$$
\begin{aligned}
\langle R_n,C_n\rangle_{L^2} &= \langle (\partial_t + L_n)C_n,C_n\rangle_{L^2}
\\& = \langle \Big[(\partial_t +L_n)\Big(\frac{|C_n|^2}{2}\Big) + \nu_N\sum_i^n |\nabla_{v_i}C_n|^2\Big], 1\rangle_{L^2}.
\end{aligned}$$
Integrating by parts in $z$, we get 
$$
\begin{aligned}
\frac12\frac{d}{dt}\| C_n\|_{L^2}^2 +\nu_N\sum_i^n\|\nabla_{v_i}C_n\|^2_{L^2}= \langle R_n,C_n\rangle_{L^2} + \frac12 \langle (\partial_t + L^*_n)f^{\otimes n}/f^{\otimes n} C_n, C_n\rangle_{L^2}.
\end{aligned}$$
Recall that the Vlasov equation $(\partial_t + \cL) f =0$, and so $(\partial_t + \cL^{\otimes n}) f^{\otimes n} =0$. It thus follows that 
\begin{equation}\label{calfn}
\begin{aligned} (\partial_t + L_n^*)f^{\otimes n}/f^{\otimes n} 
&= \frac{1}{N-1}\sum_{i\not =j}^n (K_{i,j} - K[f](x_i))\cdot \nabla_{v_i} \log f (z_i) + \nu_N \sum_{i=1}^n \frac{1}{f(z_i)}\Delta_{v_i}f(z_i)
\\&= \frac{1}{N-1}\sum_{i\not =j}^n V_{i,j}
+ \nu_N \sum_{i=1}^n B_i.
\end{aligned}\end{equation}
Finally, integrating by parts in $v_i$ and using the Cauchy--Schwarz inequality we can bound 
$$
\begin{aligned} 
\nu_N \sum_i^n \langle B_i , C_n^2\rangle_{L^2} 
&\le 2\nu_N \sum_i^n \langle \nabla_{v_i}\log f(z_i) \cdot \nabla_{v_i}C_n , C_n\rangle_{L^2}
\\
&\le \frac12\nu_N \sum_i^n \|\nabla_{v_i}C_n\|^2_{L^2} + C_0\nu_N n \| ( \nabla_{v_1} \log f(t,z_1)) C_n \|_{L^2}^2\nonumber\\
&\le  \frac12\nu_N \sum_i^n \|\nabla_{v_i}C_n\|^2_{L^2} + C_0 P_{T}^d \nu_N n N^{\epsilon}.
\end{aligned}$$
Here we also used the assumption  
$$\| (\nabla_v \log f)f^{\epsilon}\|_{L^\infty} \lesssim_{\epsilon}  \|\nabla_v (f^{\epsilon})\|_{L^\infty} \lesssim_{\epsilon,T}    1$$
together with H\"older's inequality and the a-priori estimates \eqref{L2bounds-apriori}:
\begin{align}\label{eq:Holder-logv}
& \int_{\DD^n} |C_n|^2 |f(t,z_1)|^{1-\epsilon} d z_1 d\mu_2 \cdots d\mu_n \nonumber\\
&\le  \left( \int_{\DD^n} |C_n|^2 f(t,z_1) d z_1 d\mu_2 \cdots d\mu_n \right)^{1-\epsilon} \left(  \int_{z_1 \in {\rm supp } f(t)}   d z_1  \left( \int_{\DD^{n-1}} |C_n|^2 d\mu_2 \cdots dz_n\right)^{1/\epsilon} \right)^{\epsilon}\nonumber\\
&\lesssim_{\epsilon,T} \|C_n\|^{2(1-\epsilon)}_{L^2(\DD^n)}   \|C_{n}\|_{L^\infty(\DD; L^2(\DD^{n-1}))}^{2\epsilon} \lesssim_{\epsilon,T}   N^{\epsilon}.
\end{align}
The above estimates hold for all $\epsilon \in (0,1)$. Under the assumption $\nu_N=N^{-\delta}$, we can choose $\epsilon=\delta/2$, so that $\nu_N N^{\epsilon}= \sqrt{\nu_N}$. The lemma follows, upon using the symmetry of $C_n$ to reduce $V_{i,j}$ to $V_{1,2}$ as stated. 
\end{proof}

Next, we modify the above strategy by integrating the equation \eqref{pLn} against $C_n \varphi_n$ with 
\begin{equation}\label{def-varphin} 
\varphi_n =\sum_{k\not =\ell}^n \frac{x_{k,\ell}\cdot v_{k,\ell}}{|x_{k,\ell}|^\theta},\quad \theta\in [0,1)
\end{equation}
where recall that $x_{k,\ell} = x_k-x_\ell$ and $v_{k,\ell} = v_k - v_\ell$. This leads to the following key input. 

\begin{lemma}\label{lem-keyL2}
Let $C_n$ be the solution to \eqref{pLn}, and introduce the weight function $w_n$ as in \eqref{weight}. Then, there is a positive constant $c_\theta>0$ depending on   $\theta \in [0,1)$ (that may vanish as $\theta \to 1$) so that 
 \begin{equation}\label{keyL2}
\begin{aligned}
c_\theta \| \sqrt{w_n} C_n\|_{L^2_{t,z}}^2 
&\le 
-2\langle R_n, \varphi_n C_n\rangle_{L^2_{t,z}} +\frac{n^4C_T}{N} \langle |V_{1,2}|, C_n^2 \rangle_{L^2_{t,z}} + C_T \nu_Nn^3(N/n)^{2\epsilon}
\\&\quad +C_Tn^2\Big[ \|C_n(T)\|_{L^2}^2+(1+\nu_N n + n^2 N^{-\epsilon} )\|C_n\|_{L^{2}_{t,z}}^2\Big]+ \nu_Nn^2 \sum_i^n \|\nabla_{v_i}C_n\|_{L^{2}_{t,z}}^2,
\end{aligned}
\end{equation}
where $C_T$ denotes some constant that depends on $T,\epsilon$ and the initial data $f_0$, but does not depend on $n, N$. Here $V_{1,2} = (K(x_1-x_2) - K[f](x_1))\cdot\nabla_{v_1}\log f(z_1)$. 
\end{lemma}
\begin{proof}
Recalling \eqref{calfn}, we compute 
\begin{equation}\label{L2weight}
\begin{aligned}
&\frac{d}{dt}\langle C_n^2, \varphi_n\rangle_{L^2} + 2\nu \sum_i^n \langle |\nabla_{v_i}C_n|^2, \varphi_n\rangle_{L^2} 
\\& = 2\langle R_nC_n, \varphi_n \rangle_{L^2} + \langle (\partial_t + L^*_n)(\varphi_nf^{\otimes n})/f^{\otimes n} C_n, C_n\rangle_{L^2}
\\&= 2\langle R_nC_n, \varphi_n \rangle_{L^2} + \langle(\partial_t + L_n^*)\varphi_n, C_n^2\rangle_{L^2}+  \nu_N \sum_{i=1}^n \langle B_i \varphi_n , C_n^2\rangle_{L^2}
\\&\quad + \frac{1}{N-1}\sum_{i\not =j}^n\langle V_{i,j}\varphi_n, C_n^2\rangle_{L^2} 
+2\nu_N \sum_i^n \langle \nabla_{v_i}\varphi_n\cdot \nabla_{v_i}\log f(z_i)C_n, C_n\rangle_{L^2},\end{aligned}
\end{equation}
which holds for any test function $\varphi_n$. To obtain good lower bounds for  $(\partial_t+L_n^*)\varphi_n$, we take $\varphi_n$ as in \eqref{def-varphin}. Observe that, for $\theta< 1$, by the Cauchy--Schwarz inequality we have
\begin{equation}
\begin{aligned}
\sum_i v_i \cdot \nabla_{x_i} \varphi_n =  \sum_{k\ne \ell} \left(  \frac{|v_{k,\ell}|^2}{|x_{k,\ell}|^\theta} -  \theta \frac{(x_{k,\ell}\cdot v_{k,\ell})^2}{|x_{k,\ell}|^{\theta+2}} \right) \ge  \sum_{k\ne \ell}  (1-\theta)\frac{|v_{k,\ell}|^2}{|x_{k,\ell}|^\theta} \ge  (1-\theta) w_n - n^2
\end{aligned}
\end{equation}
with the weight $w_n$ as in \eqref{weight}. Moreover, using  
$$\begin{aligned}
\frac{1}{N-1}\sum_{i\not =j}^n K_{i,j}\cdot \nabla_{v_i} \frac{x_{k,\ell}\cdot v_{k,\ell}}{|x_{k,\ell}|^\theta}  = \frac{1}{N-1} \left(\sum_{j\not =k}^n K_{k,j}  - \sum_{j\not =l}^n K_{l,j} \right)\cdot \frac{x_{k,\ell}}{|x_{k,\ell}|^\theta} 
\end{aligned}$$
for all $k\ne \ell $, and the assumption  $|K(x)|\le \frac{C_K}{|x|^{\alpha+1}}$, we get
$$\begin{aligned}
\frac{1}{N-1} \left| \sum_{i\not =j}^n K_{i,j}\cdot \nabla_{v_i} \varphi_n\right| \le C_T   \frac{n}{N} \sum_{k\ne \ell} \frac{1}{|x_{k,\ell}|^{\alpha+1}}
\end{aligned}$$
on the support of $f(t)^{\otimes n}$ with $t\in [0,T]$. Here we used that $|x| \le P_T$ in the support of $f(t)$ for $t\in [0,T]$, see \eqref{defmaxPT}. 

Similarly, we compute 
$$ \begin{aligned}
 \frac{N-n}{N-1}\sum_{i=1}^n K[f](x_i)\cdot \nabla_{v_i} \varphi_n &=  \frac{N-n}{N-1}\sum_{k\not =\ell}^n (K[f](x_k) - K[f](x_\ell))\cdot\frac{x_{k,\ell}}{|x_{k,\ell}|^\theta}  
 \\
 &=  \frac{2(N-n)}{N-1}\sum_{k\not=\ell}^n K[f](x_k) \cdot \frac{x_{k,\ell}}{|x_{k,\ell}|^\theta}  
\end{aligned}
$$
in which we recall that $K[f] = K\star\int f\; dv$. We may bound 
$$ \begin{aligned}
\Big| \frac{N-n}{N-1}\sum_{i=1}^n K[f](x_i)\cdot \nabla_{v_i} \varphi_n\Big| \le  2\|K[f]\|_{L^\infty}\sum^n_{k\le \ell}|x_{k,\ell}|^{1-\theta}  \le C_T n^2
\end{aligned}
$$
on the support of $f(t)^{\otimes n}$ with $t\in [0,T]$. 

Finally, we note that $\Delta_{v_i}\varphi_n =0$ by the construction of $\varphi_n$ in \eqref{def-varphin}. Therefore, combining the previous estimates, and recalling the weight $w_n$ as in \eqref{weight}, we obtain 
\begin{equation}\label{goodlowbd}
(\partial_t + L_n^* )\varphi_n \ge (1-\theta) w_n - C_T \frac{n}{N} \sum_{k\ne \ell} \frac{1}{|x_{k,\ell}|^{\alpha+1}} - C_T n^2
\end{equation}
on the support of $f(t)^{\otimes n}$ with $t\in [0,T]$. 

Putting this into \eqref{L2weight}, we have 
\begin{align}\label{ineqC1}
(1-\theta) \langle w_n, C_n^2 \rangle_{L^2( \DD^n)}
&\le 
\frac{d}{dt}\langle \varphi_n, C_n^2 \rangle_{L^2} - 2\langle R_n, \varphi_n C_n\rangle_{L^2} + C_T n^2 \|C_n\|_{L^2}^2  \nonumber
\\
&\quad - \frac{1}{N-1}\sum_{i\not =j}^n\langle V_{i,j}\varphi_n, C_n^2\rangle_{L^2}   + \nu_N\sum_i^n \langle |\nabla_{v_i}C_n|^2 , \varphi_n\rangle_{L^2} 
\nonumber\\& \quad  +C_T \frac{n}{N} \sum_{k\ne \ell} \langle |x_{k,\ell}|^{-\alpha-1}, C_n^2 \rangle_{L^2}     -  \nu_N \sum_{i=1}^n \langle B_i \varphi_n , C_n^2\rangle_{L^2}.
\end{align}

Let us estimate further the right-hand side of \eqref{ineqC1}. We keep the first, second and third terms. Recall that the $L^2$-integration is with respect to the probability measure $d\mu_{[n]} = f^{\otimes n} dz_{[n]}$, where $f$ is compactly supported in both $x,v$. In particular, $|x|+|v|\le P_T$ for $x,v$ in the support of $f(t)$ for $t\in [0,T]$, see \eqref{defmaxPT}. Therefore, with any choice $\theta\in [0,1]$ we have 
\begin{equation}\label{bdvarphin}
 |\varphi_n| \le \sum_{k,\ell}^n |x_{k,\ell}|^{1-\theta} |v_{k,\ell}| \le  P^2_T n^2,
 \end{equation}
 on the support of $f(t)^{\otimes n}$. Moreover, note that the most singular part of $V_{ij}$ is $|K_{ij}|\lesssim |x_{ij}|^{-\alpha-1}$ since $\|K[f]\|_{L^\infty}$ is bounded. Therefore, using also the symmetry of $C_n$, we bound the fourth and fifth terms on the  hand side of \eqref{ineqC1} by
$$
\frac{1}{N-1}\sum_{i\not =j}^n |\langle V_{i,j}\varphi_n, C_n^2\rangle_{L^2} | \le \frac{P_T^2n^4}{N-1}  \langle |V_{1,2}|, C_n^2 \rangle_{L^2} 
$$
and 
$$
 \nu_N\sum_i^n \langle |\nabla_{v_i}C_n|^2 , \varphi_n\rangle_{L^2} \le \nu_N n^2 \sum_i^n \|\nabla_{v_i}C_n\|^2_{L^2}.
$$

Next, to control $\langle |x_{k,\ell}|^{-\alpha-1}, C_n^2 \rangle_{L^2}$, we choose $\max\{0,2\alpha+2-d\}<\theta<1$.  This ensures that the function $ \left( 1+ \frac{|v|^2}{|x|^{\theta}}\right)^{-1}  |x|^{-2\alpha-2}$  belongs to $L^{1+{\epsilon}}(\DD,d\mu)$ with $\epsilon>0$ sufficiently small. Moreover, from the a-priori estimates in \eqref{L2bounds-apriori} and H\"older's inequality we can bound 
\begin{equation}\label{eq:Cn-Holder}
 \| C_n\|_{L^{(1+\epsilon)/\epsilon}(\DD^2; L^2(\DD^{n-2}))} \le  \| C_n\|_{L^2 (\DD^{n})}  \| C_n\|_{L^\infty(\DD^2; L^2(\DD^{n-2}))}^{1/(1+\epsilon)} \le N^{1/(1+\epsilon)} \| C_n\|_{L^2 (\DD^{n})} . 
\end{equation}
Then by using the symmetry of $C_n$ and the Cauchy--Schwarz inequality we have
\begin{equation*}
\begin{aligned}
\frac{n}{N}\langle |x_{k,\ell}|^{-\alpha-1}, C_n^2 \rangle_{L^2}  &\le \frac{1}{N} \left \langle n^2\left( 1+ \frac{|v_{k,\ell}|^2}{|x_{k,\ell}|^{\theta}}\right), C_n^2 \right\rangle_{L^2}^{1/2}   \left \langle \left( 1+ \frac{|v_{k,\ell}|^2}{|x_{k,\ell}|^{\theta}}\right)^{-1}  |x_{k,\ell}|^{-2\alpha-2}, C_n^2 \right\rangle_{L^2}^{1/2}\\
&\le \frac{C}{N} \| \sqrt{w_n} C_n\|_{L^2} \left\| \left( 1+ \frac{|v_{k,\ell}|^2}{|x_{k,\ell}|^{\theta}}\right)^{-1}  |x_{k,\ell}|^{-2\alpha-2} \right\|^{1/(1+\epsilon)}_{L^{1+\epsilon}(\DD^2)}  \| C_n\|_{L^{(1+\epsilon)/\epsilon}(\DD^2; L^2(\DD^{n-2}))} \\
&\le \frac{C}{N^{\epsilon/(1+\epsilon)}}   \| \sqrt{w_n} C_n\|_{L^2}  \| C_n\|_{L^2 (\DD^{n})} 
\end{aligned}
\end{equation*}
for all $k\ne \ell$. Thus, by the Cauchy--Schwarz inequality again,
\begin{align}\label{eq:singular-CS}
C_T\frac{n}{N}\sum_{k\ne \ell}^n\langle |x_{k,\ell}|^{-\alpha-1}, C_n^2 \rangle_{L^2} 
 &\le  \frac{C_T n^2}{N^{\epsilon/(1+\epsilon)}}   \| \sqrt{w_n} C_n\|_{L^2}  \| C_n\|_{L^2 (\DD^{n})}  \nonumber \\
 &\le \frac{1-\theta}{2}   \| \sqrt{w_n} C_n\|_{L^2}^2 + \frac{C_T^2 n^4}{N^{2\epsilon/(1+\epsilon)}} \| C_n\|_{L^2 (\DD^{n})} ^2.
\end{align}
The term $\frac{1-\theta}{2}   \| \sqrt{w_n} C_n\|_{L^2}^2 $ in \eqref{eq:singular-CS}  will be cancelled by (half of) the left-hand side of    
\eqref{ineqC1}. 

For the last term on the right-hand side of \eqref{ineqC1}, recalling that
$B_i=f(z_i)^{-1}\Delta_{v_i}f(z_i)$. By integration
by parts in $v_i$ we can write 
\[
\begin{aligned}
\nu_N\langle B_i\varphi_n,C_n^2\rangle_{L^2}
={}&-\nu_N\langle
\nabla_{v_i}\log f(z_i)\cdot\nabla_{v_i}\varphi_n,
C_n^2\rangle_{L^2}\\
&-2\nu_N\langle
\nabla_{v_i}\log f(z_i)\cdot\nabla_{v_i}C_n\,
\varphi_n,C_n\rangle_{L^2}.
\end{aligned}
\]
We can use  $|\varphi_n|\le C_Tn^2$ and
$|\nabla_{v_i}\varphi_n|\le C_Tn$ in the support of $f^{\otimes n}$. Moreover, although $\nabla_v\log f$ is not uniformly bounded, we have the assumption $
\|f^\epsilon\nabla_v\log f\|_{L^\infty}
\lesssim_{\epsilon,T} 1$. Consequently, by Hölder's inequality and
the a priori bounds \eqref{eq:Cn-Holder},
\[
\begin{aligned}
\|(\nabla_{v_i}\log f(z_i))C_n\|_{L^2}^2
&\lesssim_{\epsilon,T}
\int_{D^n} C_n^2 f(z_i)^{-2\epsilon}
\,d\mu_{[n]}\\
&\lesssim_{\epsilon,T} 
\|C_n\|_{L^2}^{2(1-2\epsilon)}
\|C_n\|_{L^\infty(\DD;L^2(\DD^{n-1}))}^{4\epsilon} \lesssim_{\epsilon,T} (N/n)^{2\epsilon}.
\end{aligned}
\]
Using $|\varphi_n|\le C_Tn^2$ and
$|\nabla_{v_i}\varphi_n|\le C_Tn$, we therefore obtain
\[
\begin{aligned}
\nu_N\sum_{i=1}^n
|\langle B_i\varphi_n,C_n^2\rangle_{L^2}|
\lesssim_{\epsilon,T}&
\nu_Nn^2
\sum_{i=1}^n\|\nabla_{v_i}C_n\|_{L^2}^2+\nu_Nn^3(N/n)^{2\epsilon}.
\end{aligned}
\]

Finally, inserting the above estimate into  \eqref{ineqC1}, then integrating in time and recalling that $C_n(0)=0$, we obtain the lemma. 
\end{proof}

\begin{proof}[Proof of Proposition \ref{prop-keyL2}] We first claim that 
\begin{equation}\label{claimL2V12}
\begin{aligned}
\frac{n^2}{N}\|\sqrt{|V_{1,2}|}C_n\|_{L^2(\DD^n)}^2 &\le C_T (N/n)^{-\delta_\alpha } \|\sqrt{w_n}C_n\|_{L^2(\DD^n)} 
\end{aligned}
\end{equation}
for some $\delta_\alpha >0$. Indeed, by the H\"older inequality, we bound 
\begin{equation}\label{L2V12}
\begin{aligned}
 \Big\|\sqrt{|V_{1,2}|}C_n\Big\|_{L^2(\DD^n)}^2
&\le 
\Big\| ( 1+ \frac{|v_{1,2}|}{|x_{1,2}|^{\theta/2}}) C_n\Big\|_{L^{2}(\DD^n)}\Big \| (1+ \frac{|v_{1,2}|}{|x_{1,2}|^{\theta/2} })^{-1}V_{1,2}C_n\Big\|_{L^2(\DD^n)} 
\end{aligned}
\end{equation}
in which we recall the notation $x_{1,2} = x_1-x_2$ and $v_{1,2} = v_1-v_2$. By symmetry of $C_n$, we have  
\begin{equation}\label{L2V12w}
\sqrt{n(n-1)}\Big\| (1+ \frac{|v_{1,2}|}{|x_{1,2}|^{\theta/2}}) C_n\Big\|_{L^{2}(\DD^n)} \le  \Big\| \sqrt{\sum_{k\not =\ell}^n (1+ \frac{|v_{k,\ell}|^2}{|x_{k,\ell}|^{\theta}}) }C_n\Big\|_{L^{2}(\DD^n)} \le \|\sqrt{w_n}C_n\|_{L^2(\DD^n)}.
\end{equation}
On the other hand, using the H\"older inequality, we bound 
$$ 
\Big \| (1+  \frac{|v_{1,2}|^2}{|x_{1,2}|^{\theta}})^{-1/2}V_{1,2}C_n\Big\|_{L^2(\DD^n)} \le \Big \|  (1+  \frac{|v_{1,2}|^2}{|x_{1,2}|^{\theta}})^{-1/2} V_{1,2}\Big\|_{L^p(\DD^2)} \| C_n\|_{L^{\frac{2p}{p-2}}(\DD^2; L^2(\DD^{n-2}))}  
$$
with $p>2$ and $1/2 = 1/p + 1/r$. Recalling that $V_{1,2} = (K(x_1-x_2) - K[f](x_1))\cdot\nabla_{v_1}\log f(z_1)$, with $|K(x)|\lesssim |x|^{-1-\alpha}$. Arguing similarly as in \eqref{intro-bdVstar}, we therefore have 
\begin{equation}\label{VLp}
\Big\| \Big(1+ \frac{|v_{1,2}|^2} {|x_{1,2}|^{\theta}}\Big)^{-1/2}V_{1,2}\Big\|_{L^p(\DD^2)} \lesssim_{\epsilon,T} 1
\end{equation}
for some $p>2$ (which is sufficiently close to 2, depending on $\alpha$ and $\epsilon$). Combining with H\"older's inequality similarly to \eqref{eq:Cn-Holder}, we obtain
\begin{equation}\label{decayV12}
\begin{aligned}
\frac{n^2}{N}\Big \|  \Big(1+  \frac{|v_{1,2}|^2}{|x_{1,2}|^{\theta}}\Big)^{-1/2} V_{1,2}C_n\Big\|_{L^2(\DD^n)} 
&\lesssim_{\epsilon,T} \frac{n^2}{N} (N/n)^{2/p}.
\end{aligned}
\end{equation}
This, together with \eqref{L2V12}-\eqref{L2V12w}, proves the claim \eqref{claimL2V12}, with $\delta_\alpha = 1-2/p$ (recalling that the choice of $p$ depends on $\alpha$ and $\epsilon$). 

We are now ready to prove the proposition, using Lemmas \ref{lem-L2basic} and \ref{lem-keyL2}. Indeed, putting the decay estimate \eqref{claimL2V12} into \eqref{L2basic} and \eqref{keyL2}, respectively, we obtain 
 \begin{equation}\label{L2basic1}
\begin{aligned}
\frac{d}{dt}\|C_n\|_{L^{2}}^2 +\nu_N\sum_i^n\|\nabla_{v_i}C_n\|^2_{L^{2}}
&\le 2\langle R_n, C_n\rangle_{L^{2}} + C_T\nu_N n 
\\&\quad + C_T (N/n)^{-\delta_\alpha } \|\sqrt{w_n}C_n\|_{L^2(\DD^n)},
\end{aligned}
\end{equation}
and
 \begin{equation}\label{keyL21}
\begin{aligned}
c_\theta \| \sqrt{w_n} C_n\|_{L^2_{t,z}}^2 
&\le 
-2\langle R_n, \varphi_n C_n\rangle_{L^2_{t,z}} +C_Tn^2(N/n)^{-\delta_\alpha } \|\sqrt{w_n}C_n\|_{L^2_{t,z}} + C_T \nu_Nn^3(N/n)^{2\epsilon}
\\&\quad +C_Tn^2\Big[ \|C_n(T)\|_{L^2}^2+(1+\nu_N n + n^2 N^{-\epsilon} )\|C_n\|_{L^{2}_{t,z}}^2\Big]+ \nu_Nn^2 \sum_i^n \|\nabla_{v_i}C_n\|_{L^{2}_{t,z}}^2,
\end{aligned}
\end{equation}
in which $C_T$ depends on $T,\alpha,\epsilon$ and $f$. Integrating \eqref{L2basic1} in time and using the Gronwall inequality, we obtain  
\begin{equation}\label{supL2Cn}
\begin{aligned}\sup_{[0,T]}\|C_n(t)\|_{L^{2}}^2 + \nu_N\sum_i^n\|\nabla_{v_i}C_n\|^2_{L^{2}_{t,z}}
&\le C_T \int_0^T |\langle R_n,C_n\rangle_{L^{2}_{z}}| dt  + C_T  (N/n)^{-\delta_\alpha } \|\sqrt{w_n}C_n\|_{L^2_{t,z}} + C_T \nu_N n
\end{aligned}\end{equation}
provided $\nu_N n\lesssim 1$. 
Putting this into \eqref{keyL21} and recalling \eqref{bdvarphin}, we obtain 
\begin{align}\label{keyL22}
c_\theta \| \sqrt{w_n} C_n\|_{L^2_{t,z}}^2 
&\le 
C_T n^2 \int_0^T |\langle R_n,C_n\rangle_{L^{2}_{z}}| dt  +
C_Tn^2(N/n)^{-\delta_\alpha } \|\sqrt{w_n}C_n\|_{L^2_{t,z}} + C_T \nu_Nn^3(N/n)^{2\epsilon} \nonumber\\
&+ C_Tn^2\Big[ \|C_n(T)\|_{L^2}^2+(1+\nu_N n + n^2 N^{-\epsilon} )\|C_n\|_{L^{2}_{t,z}}^2\Big].
\end{align}
On the other hand, we can simply bound $|\langle R_n, C_n\rangle_{L^2_{z}}| \le  \| C_n\|_{L^2_{z}}\|R_n\|_{L^{2}_{z}}.
$
Therefore, combining \eqref{supL2Cn} and \eqref{keyL22} and introducing 
$$\zeta_T = \Big[ n^2 \sup_{[0,T]}\|C_n(t)\|_{L^{2}}^2 + n^2\nu_N\sum_i^n\|\nabla_{v_i}C_n\|^2_{L^{2}_{t,z}} + 
c_\theta \| \sqrt{w_n} C_n\|_{L^2_{t,z}}^2  \Big]^{1/2} $$
we deduce 
$$\zeta_T^2 \le C_T n  \|R_n\|_{L^{2}_{t,z}}\zeta_T + C_Tn^2(N/n)^{-\delta_\alpha } \zeta_T + C_T \nu_Nn^3(N/n)^{2\epsilon}.$$
With $\nu_N=N^{-\delta}$, this yields the boundedness of $\zeta_T$, and thus Proposition \ref{prop-keyL2}. 
\end{proof}

\subsection{$H^{-1}$ sources}\label{sec:H-1sources}

In this section, we derive weighted estimates on the solution $C_n$ to \eqref{pLn} with $H^{-1}$ sources. Specifically, we obtain the following proposition. 

\begin{proposition}\label{prop-keydvS} Let $K(x)$ be the Riesz interaction force \eqref{Riesz} with $0\le \alpha < \frac{d-1}{2}$, and let $C_n$ be the solution to \eqref{pLn} with 
\begin{equation}\label{RdvS}
R_n = \sum_{i=1}^n \nabla_{v_i}\cdot S_{n,i},
\end{equation}
 and satisfy the a priori bounds in Lemma \ref{lem-apriori}. Then, for $\nu_N n\lesssim 1$, there are positive constants $\delta_\alpha$ and $C_{\alpha,T}$, independent of $n$, so that 
\begin{equation}
\label{keysupL2}
\begin{aligned}
n^2\sup_{[0,T]}&\|C_n(t)\|_{L^{2}(\DD^n)}^2+ \nu_Nn^2\sum_i^n\|\nabla_{v_i}C_n\|^2_{L^{2}([0,T]\times \DD^n)} +\| \sqrt{w_n} C_n\|_{L^2([0,T]\times \DD^n)}^2 
\\&\le C_{\alpha,T} \Big(n^2\nu^{-1}\sum_i^n\| S_{n,i}\|_{L^2([0,T]\times \DD^n)}^2 + n^4(N/n)^{-2\delta_\alpha }\Big),
\end{aligned}
\end{equation}
for all $1\le n\le N$, where the weight $w_n$ is defined as in \eqref{weight}. 
\end{proposition}
\begin{proof} As in the proof of Proposition \ref{prop-keyL2}, we start from the two estimates \eqref{supL2Cn}-\eqref{keyL22}. Using \eqref{RdvS}, we 
 compute the integration by parts
\begin{equation}\label{eq:bddvS-0}
 \langle R_n,C_n\rangle_{L^{2}_z} =  -\sum_{i=1}^n \langle S_{n,i}, \nabla_{v_i}C_n\rangle_{L^{2}_z}  - \sum_{i=1}^n  \langle S_{n,i}, C_n \nabla_{v_i} \log f (z_i) \rangle_{L^{2}_z} .
 \end{equation} 
Applying the Cauchy--Schwarz inequality for the first term on the left-hand side of \eqref{eq:bddvS-0}, and also using the assumption $f^{\epsilon} |\nabla_v\log f|\lesssim 1$, 
and H\"older's inequality as in \eqref{eq:Holder-logv} for the second term, we find that 
\begin{equation}\label{bddvS}
\begin{aligned} 
|\langle R_n,C_n\rangle_{L^{2}_z}| \le    \Big(\sum_i\| \nabla_{v_i}C_n\|_{L^{2}}^2 \Big)^{1/2}\Big(\sum_i\| S_{n,i}\|_{L^{2}}^2\Big)^{1/2} + C_{T} \sqrt{n} (N/n)^{\epsilon} \Big(\sum_i\| S_{n,i}\|_{L^{2}}^2\Big)^{1/2}.
\end{aligned}\end{equation}
Therefore, putting this into \eqref{supL2Cn}, we obtain 
\begin{equation}\label{supL2Cn1}
\begin{aligned}
\sup_{[0,T]}\|C_n(t)\|_{L^{2}}^2 + \nu_N\sum_i^n\|\nabla_{v_i}C_n\|^2_{L^{2}_{t,z}}
&\le \Big(\sum_i\| \nabla_{v_i}C_n\|_{L^{2}_{t,z}}^2 \Big)^{1/2}\Big(\sum_i\| S_{n,i}\|_{L^{2}_{t,z}}^2\Big)^{1/2} 
\\&\quad + C_T (N/n)^{-\delta_\alpha } \|\sqrt{w_n}C_n\|_{L^2_{t,z}} + C_T \nu_N n
\end{aligned}\end{equation}
provided that $\nu_Nn \lesssim 1$. 
Similarly, from \eqref{keyL22}, we obtain 
\begin{align}\label{keyL23}
c_\theta \| \sqrt{w_n} C_n\|_{L^2_{t,z}}^2 
&\le C_0 n^2 \Big(\sum_i\| \nabla_{v_i}C_n\|_{L^{2}_{t,z}}^2 \Big)^{1/2}\Big(\sum_i\| S_{n,i}\|_{L^{2}_{t,z}}^2\Big)^{1/2} 
+
C_Tn^2(N/n)^{-\delta_\alpha } \|\sqrt{w_n}C_n\|_{L^2_{t,z}} \nonumber\\
&+ C_T \nu_Nn^3(N/n)^{2\epsilon} + C_Tn^2\Big[ \|C_n(T)\|_{L^2}^2+(1+\nu_N n + n^2 N^{-\epsilon} )\|C_n\|_{L^{2}_{t,z}}^2\Big].
\end{align}
Next, introducing 
$$\zeta_T = \Big[ C_Tn^2 \sup_{[0,T]}\|C_n(t)\|_{L^{2}}^2 + n^2\nu_N\sum_i^n\|\nabla_{v_i}C_n\|^2_{L^{2}_{t,z}} + 
c_\theta \| \sqrt{w_n} C_n\|_{L^2_{t,z}}^2 \Big]^{1/2} $$
and, from combining \eqref{supL2Cn1} and \eqref{keyL23}, deducing  
$$\zeta_T^2 \le C_T n \nu^{-1/2}\Big(\sum_i\| S_{n,i}\|_{L^{2}_{t,z}}^2\Big)^{1/2}  \zeta_T + C_Tn^2(N/n)^{-\delta_\alpha } \zeta_T + C_T \nu_Nn^3(N/n)^{2\epsilon}$$
which yields the boundedness of $\zeta_T$, 
we obtain the proposition.  
\end{proof}

\section{Iterative bounds}\label{sec:Iterative bounds}

In this section, we return to bound the adjoint cumulants $C_{N,n}$. We denote $C_{N,n}=C_n$ for simplicity (dropping the apparent dependence on $N$) whose dynamics is 
obtained in Proposition \ref{prop-CNn}. Specifically, we obtain the following proposition. 

\begin{proposition}\label{prop-main} Let $d\ge 2$. 
Let $C_n=C_{N,n}$ be the adjoint cumulants as constructed in Proposition \ref{prop-CNn}. Then, for all $n\ge 1$, there hold 
\begin{align}\label{keyfullCn}
\langle (\partial_t + L_n)C_n, C_n \rangle_{L^2_z}& \le C_T\Big[\| g_n\|_{L^2}+n (N/n)^{-1/2 + 1/d} + \nu_N \sqrt{Nn}\Big] \| C_n\|_{L^2}\nonumber\\
&\quad + \frac{C_T}{\sqrt{N}} \Big( \| \sqrt{w_{n}} C_{n}\|_{L^2(\DD^{n})} + \| \sqrt{w_{n+1}} C_{n+1}\|_{L^2(\DD^{n})} \Big) \times \nonumber\\
&\qquad \times \Big[ \Big( \sum_i^n \|\nabla_{v_i}C_n\|_{L^2}^2\Big)^{1/2} + \Big( \sum_{i}^n \|\nabla_{v_i}C_{n-1}\|^2_{L^2(\DD^{n-1})} \Big)^{1/2} +  \sqrt{n}(N/n)^{\epsilon}\Big]\nonumber\\
&\quad + C_T \Big(n\|C_n\|_{L^2(\DD^{n})}
+ \| \sqrt{w_{n+2}}C_{n+2}\|_{L^2(\DD^{n+2})}\Big) \|C_n\|_{L^2},
\end{align}
for some constant $C_T$ that is independent of $n$. In addition, letting $\varphi_n =\sum_{k\not =\ell}^n \frac{x_{k,\ell}\cdot v_{k,\ell}}{|x_{k,\ell}|^\theta}$, we obtain the same bounds for $\langle (\partial_t + L_n)C_n, \varphi_n C_n \rangle_{L^2}$, upon replacing $C_T$ in \eqref{keyfullCn} by $C_T n^2$. 
\end{proposition}

Proposition \ref{prop-main} provides bounds on $C_n$ that solves the full cumulant dynamics. Observe that the last two terms on the right are of order one in the large $N$, one of which involves higher-order cumulants $C_{n+2}$, not to mention the singular weight $w_n$. The crucial observation is that the estimate \eqref{keyfullCn} in fact holds for weighted norms, up to some loss in $n$, which play a crucial role in closing the iteration. 

\subsection{Cumulant dynamics}

In view of Proposition \ref{prop-CNn} and the linear analysis established in Section \ref{sec-linear}, we study the structure of the cumulant dynamics. Precisely, we shall write the dynamics of cumulants $C_n$ as 
\begin{equation}\label{fullCn} (\partial_t + L_n)C_n = R_n
\end{equation}
with zero initial data $C_n(0)=0$, with the same linear operator $L_n$ as in Proposition \ref{prop-CNn}, and  
\begin{equation}\label{Rnstructure}
\begin{aligned}
R_n &= g_{n}  + \TT_{n,-}C_{n-1} + \TT_{n,0}C_{n} + \TT_{n,1}C_{n+1} + T_{n,2}C_{n+2} + D_{n,-}C_{n-1} + D_{n,0}C_{n} + D_{n,1}C_{n+1} 
\end{aligned}
\end{equation}
in which $D_{n,*}$ collect all the terms in $T_{n,*}$ that are of a divergence form, and $\TT_{n,*}$ denote the remaining terms, namely $\TT_{n,*} = T_{n,*} - D_{n,*}$. Precisely, in view of Proposition \ref{prop-CNn}, we have set 
\begin{equation}\label{defDnstar}
\begin{aligned}
D_{n,-}C_{n-1} & = a_{n,-}\sum_{i\not = j}^n \nabla_{v_i} \cdot [(K[f](x_i) - K_{i,j}))C_{n-1}(z_{[n]\setminus\{j\}})] 
\\
D_{n,0} C_{n} & = 
\frac{1}{N-1}\sum_{i\not =j}^n \nabla_{v_i}\cdot \int_{\DD} K_{i,*} C_n(z_{[n]\setminus\{j\}}, z_*)\; d\mu_*
\\
D_{n,1} C_{n+1} & = 
- a_{n,1}
\sum_{i=1}^n 
\nabla_{v_i}\cdot \int_{\DD} K_{i,*}C_{n+1}(z_{[n]}, z_*, )
\; d\mu_*
\end{aligned}
\end{equation}
with $a_{n,-}\le \frac{2}{\sqrt{nN}}$ and $a_{n,1}\le 2\sqrt{n/N}$. Note that there are no derivative terms in $\TT_{n,*}$ and in $T_{n,2}C_{n+2}$. We shall apply the results obtained in Propositions \ref{prop-keyL2} and \ref{prop-keydvS}. This reduces to derive estimates on $R_n$ in $L^2$ or $H^{-1}$ norms. 

\subsection{$L^2$ components}

In this section, we bound the source terms that have no derivatives in $R_n$. Namely, we have the following. 

\begin{lemma}\label{lem-Cnzero} Let $\TT_{n,*}$ be defined as in \eqref{Rnstructure}, and set
$$R_n =\TT_{n,-}C_{n-1} + \TT_{n,0}C_{n} + \TT_{n,1}C_{n+1}$$
with $C_n$ satisfying the a priori bounds in Lemma \ref{lem-apriori}. Then, there are positive constants $\delta_\alpha$ and $C_{\alpha,T}$ independent of $n$ so that
\begin{equation}
\label{bdCnzero}
\begin{aligned}
\| R_n - T_{n,0,1}C_n\|_{L^2(\DD^n)} \le C_Tn (N/n)^{-1/2 + 1/d} + C_T \nu_N \sqrt{Nn},
\end{aligned}
\end{equation}
in which 
\begin{equation}\label{orderone}  T_{n,0,1}C_n = \frac{N-n}{N-1} \sum_{j=1}^n \int_{\DD} V_{*,j}C_n(z_{[n]\setminus\{j\}}, z_*)\; d\mu_* .
\end{equation}
\end{lemma}

\begin{proof} Recall that $\TT_{n,*} = T_{n,*} - D_{n,*}$, namely that of $T_{n,*}$ defined as in Proposition \ref{prop-CNn} without the derivative terms. For convenience, we recall    
$$\begin{aligned}
\TT_{n,0} C_{n} & = \frac{N-n}{N-1} \sum_{j=1}^n \int_{\DD} V_{*,j}C_n(z_{[n]\setminus\{j\}}, z_*)\; d\mu_* 
- \frac{1}{N-1}\sum_{i\not =j}^n  \int_{\DD} V_{*,j} C_n(z_{[n]\setminus\{i\}}, z_*)\; d\mu_*
\\&\quad + \frac{1}{N-1}\sum_{i\not =j}^n  \int_{\DD^2} K_{*,*'} C_n(z_{[n]\setminus\{i,j\}}, z_*, z_*')\; d\mu_*d\mu_*'
\\
\TT_{n,-}C_{n-1} & = 
- a_{n,-} \sum_{i\not =j}^n \int_{\DD} V_{*,j} C_{n-1}(z_{[n]\setminus\{i,j\}}, z_*)\; d\mu_{*}
\\
\TT_{n,1} C_{n+1} & = 
a_{n,1}
\sum_{j=1}^n 
\int_{\DD} V_{*,j}C_{n+1}(z_{[n]}, z_* )
\; d\mu_*  -2a_{n,1} \sum_{i=1}^n\int_{\DD^2}V_{*,*'} C_{n+1}(z_{[n]\setminus\{i\}}, z_*, z_*') d\mu_* d\mu_*'
\\&\quad  + \nu_{n,1} \int_{\DD} B_*C_{n+1}(z_{[n]}, z_*)
\; d\mu_*
\end{aligned}
$$
with $a_{n,-}= \frac{\sqrt{N-n+1}}{\sqrt{n}(N-1)}$, $a_{n,1}= \frac{\sqrt{(N-n)(n+1)}}{N-1}$, and $\nu_{n,1} = \nu_N \sqrt{(N-n)(n+1)}$. In what follows, we denote by $\TT_{n,*,j}$, $j=1,2,3$, the respective terms in $\TT_{n,*}$. In particular, the first term in $\TT_{n,0}$ defines $T_{n,0,1}$ as in \eqref{orderone}, which is of order one and will be treated in the next section. Let us bound the remaining terms. Indeed, 
we first bound $\TT_{n,0,2}$, namely the second integral term in $\TT_{n,0}$. Note that the local singularity of $V_{*,j}$ is of order $|x_{*,j}|^{-1-\alpha}$, and so $V_{*,j} \in L^r(\DD_*)$ with $1\le r < p_* = \frac{d}{1+\alpha}$. As $\alpha <\frac{d-1}{2}$, we have $p_* > \frac{2d}{d+1}$. In particular, $V_{*,j} \in L^{\frac{2d}{d+1}}(\DD_*)$. Therefore, we bound 
$$
\begin{aligned}
|\TT_{n,0,2}C_n| 
&\le  \frac{1}{N-1}\sum_{j\not = i}^n\int_{\DD} V_{*,j} C_n(z_{[n]\setminus\{i\}}, z_*)\; d\mu_* 
\\& \le  \frac{1}{N-1}\sum_{j\not = i}^n\| V_{*,j}^{\frac{d}{d+1}}\|_{L^2(\DD_*)}\| V_{*,j}^{\frac{1}{d+1}} C_n (z_{[n]\setminus\{i\}}, z_*)\|_{L^2(\DD_*)} 
\\&\le 
 \frac{C}{N-1}\sum_{j\not = i}^n\| V_{*,j}^{\frac{1}{d+1}} C_n(z_{[n]\setminus\{i\}}, z_*)\|_{L^2(\DD_*)}.
\end{aligned}$$
Hence, by triangle inequality, we bound 
$$ 
\begin{aligned}
\| \TT_{n,0,2}C_n\|_{L^2(\DD^n)}
&\le \frac{1}{N-1}\sum_{j\not = i}^n\Big\|  V_{*,j}^{\frac{1}{d+1}} C_n(z_{[n]\setminus\{i\}}, z_*)
\Big\|_{L^2([0,T]\times \DD^{n+1})}
\\
&\le \frac{1}{N-1}\sum_{j\not = i}^n\|V_{*,j}\|^{\frac{1}{d+1}}_{L^2_TL^{\frac{2d}{d+1}}(\DD^2_{*,j})}  \|C_n\|_{L^2_T L^{\frac{2d}{d-1}}(\DD^2_{*,j}; L^2(\DD^{n-2}))}
\\
&\le \frac{n^2}{N-1}\|C_n\|^{1/d}_{L^2_T L^{\infty}(\DD^2_{*,j}; L^2(\DD^{n-2}))}\|C_n\|^{1-1/d}_{L^2_T L^2(\DD^{n})}
.\end{aligned}
$$
Therefore, using the a priori bounds from Lemma \ref{lem-apriori}, we have 
\begin{equation}\label{bdTTn02}\| \TT_{n,0,2}C_n\|_{L^2(\DD^n)} \le C_Tn (N/n)^{-1 + 1/d}.
\end{equation} 
Similarly, since $K_{*,*'} \in L^{\frac{2d}{d+1}}(\DD^2_{*,*'})$, we bound 
$$
\begin{aligned}
|\TT_{n,0,3}C_n| &\le \frac{1}{N-1}\sum_{i\not =j}^n  \int_{\DD^2} K_{*,*'} C_n(z_{[n]\setminus\{i,j\}}, z_*, z_*')\; d\mu_*d\mu_*' 
\\&\le \frac{C}{N-1}\sum_{j\not = i}^n\| K_{*,*'}^{\frac{1}{d+1}} C_n(z_{[n]\setminus\{i,j\}}, z_*,z_{*'})\|_{L^2(\DD_{*,*'}^2)},
\end{aligned}$$
and so, $\| \TT_{n,0,3}C_n\|_{L^2} \le C_Tn (N/n)^{-1 + 1/d}$, similarly as was done for $\TT_{n,0,2}C_n$.  

Next, we observe that $\TT_{n,-}C_{n-1}$ and $\TT_{n,1,1}C_{n+1}$ are of a similar form to $\TT_{n,0,2}C_n$, while $\TT_{n,1,2} C_{n+1}$ is of a similar form to $\TT_{n,0,3}C_n$. Hence, following the above analysis and noting that $a_{n,-}\le \frac{2}{\sqrt{nN}}$ and $a_{n,1}\le 2\sqrt{n/N}$, we have 
\begin{equation}\label{bdTTn123}
\begin{aligned}
\| \TT_{n,-}C_{n-1}\|_{L^2(\DD^n)} + \| \TT_{n,1,1}C_{n+1}\|_{L^2(\DD^n)}+ \| \TT_{n,1,2}C_{n+1}\|_{L^2(\DD^n)} &\le C_Tn(N/n)^{-1/2 + 1/d}.
\end{aligned}
\end{equation}
Finally, concerning the integral term $\TT_{n,1,3}C_{n+1}$, we bound 
$$\begin{aligned} \|\TT_{n,1,3}C_{n+1}\|_{L^2(\DD^n)} 
&\le  \nu_N \sqrt{(N-n)(n+1)}\int_{\DD} |B_*| \|C_{n+1}(z_{[n]}, z_*, )\|_{L^2(\DD^n)}
\; d\mu_* 
\\&\le C_T \nu_N \sqrt{Nn}
\| C_{n+1}\|_{L^2(\DD^{n+1})},
\end{aligned}$$
which is bounded by $C_T \nu_N \sqrt{Nn}$, upon using the a priori bounds from Lemma \ref{lem-apriori}. The lemma thus follows. 
\end{proof}

\subsection{$H^{-1}$ components}

In this section, we bound sources $R_n$ that are in $H^{-1}$. Namely, we obtain the following. 

\begin{lemma}\label{lem-Cnminus} Let $D_{n,*}$ be defined as in \eqref{defDnstar}, and set 
$$R_n =D_{n,-}C_{n-1} + D_{n,0}C_{n} + D_{n,1}C_{n+1}$$
with $C_n$ satisfying the a priori bounds in Lemma \ref{lem-apriori}. 
Then, for each $n\ge 1$, there hold
\begin{equation}
\label{bdCnminus}
\begin{aligned}
|\langle R_n, C_n \rangle_{L^2_z}| 
&\le \frac{C_T}{\sqrt{N}} \Big( \| \sqrt{w_{n}} C_{n}\|_{L^2(\DD^{n})} + \| \sqrt{w_{n+1}} C_{n+1}\|_{L^2(\DD^{n})} \Big) \times \\
&\quad \times \Big[ \Big( \sum_i^n \|\nabla_{v_i}C_n\|_{L^2}^2\Big)^{1/2} + \Big( \sum_{i}^n \|\nabla_{v_i}C_{n-1}\|^2_{L^2(\DD^{n-1})} \Big)^{1/2} +  \sqrt{n}(N/n)^{\epsilon}\Big]
\end{aligned}
\end{equation}
for some constant $C_T$ that is independent of $n$. Alternatively, for each $n\ge 1$ we also have 
\begin{equation}
\label{bdCnminus1}
\begin{aligned}
|\langle R_n, C_n \rangle_{L^2_z}| 
&\le C_T n^{1/2}(N/n)^{-1/2+1/d} \times \\
&\quad\times\Big[ \Big( \sum_i^n \|\nabla_{v_i}C_n\|_{L^2}^2\Big)^{1/2} +\Big( \sum_{i}^n \|\nabla_{v_i}C_{n-1}\|^2_{L^2(\DD^{n-1})} \Big)^{1/2} + \sqrt{n}(N/n)^{\epsilon} \Big]. 
\end{aligned}
\end{equation}
\end{lemma}

Note that  the factor $(N/n)^{-1/2+1/d}$ does not give the smallness when $d=2$. However, this will be compensated by taking the summation of large $n$ in Section \ref{sec:final-cumulant}.

\begin{proof} In view of \eqref{defDnstar}, we write $D_{n,0}C_n + D_{n,1}C_{n+1} = \sum_{i=1}^n \nabla_{v_i}\cdot S_{n,i}$ for each respective integral terms $S_{n,i}$ of $D_{n,0}$ and $D_{n,1}$. Precisely, we have 
\begin{equation}\label{defSnstar}
\begin{aligned}
S_{n,0,i} C_{n} & = 
\frac{1}{N-1}\sum_{j = 1, \; j \not =i}^n \int_{\DD} K_{i,*} C_n(z_{[n]\setminus\{j\}}, z_*)\; d\mu_*
\\
S_{n,1,i} C_{n+1} & = 
- a_{n,1} \int_{\DD} K_{i,*}C_{n+1}(z_{[n]}, z_* )
\; d\mu_*,
\end{aligned}
\end{equation}
with $a_{n,1}=\frac{\sqrt{(N-n)(n+1)}}{N-1}$. Let us bound each term. Indeed, similar to the derivation of \eqref{bdTTn123}, we have  
\begin{equation}\label{bdSprior}
\| S_{n,0,i} C_{n} \|_{L^2_z}\le C_T (N/n)^{-1+1/d}, \qquad \| S_{n,1,i} C_{n+1} \|_{L^2_z}\le C_T (N/n)^{-1/2+1/d}.
\end{equation}
However, such an estimate is not sufficient, since we need a faster decay in $N$ in order to compensate with the vanishing of $\nu_N$. Alternatively, we bound 
\begin{equation}\label{bdKoneint} \begin{aligned}
\int_{\DD} K_{i,*}C_{n+1}(z_{[n]}, z_*)
\; d\mu_*
& \le 
\Big(
\int_{\DD} (1+ \frac {|v_{i,*}|^2}{|x_{i,*}|^\theta})^{-1}|K_{i,*}|^2\; d\mu_*\Big)^{1/2} \Big(  \int_{\DD}
(1+\frac{|v_{i,*}|^2}{|x_{i,*}|^\theta} ) |C_{n+1}(z_{[n]}, z_*)|^2
\; d\mu_*\Big )^{1/2}.
\end{aligned}
 \end{equation}
Recalling \eqref{Riesz}, we have $K_{i,*} \sim |x_{i,*}|^{-\alpha -1}$, and hence 
\begin{equation}\label{bdVstar}\int_{\DD} (1+ \frac {|v_{i,*}|^2}{|x_{i,*}|^\theta})^{-1}|K_{i,*}|^2\; d\mu_* \le C_T,
\end{equation}
provided that $2\alpha + 2 - \theta < d$. The latter condition is valid if we choose $\max\{0,2\alpha+2-d\}<\theta<1$, which is possible
under the assumption $\alpha<(d-1)/2$. Therefore, using the symmetry of $C_{n+2}$ and recalling the weight $w_n$ defined as in \eqref{weight}, we bound
\begin{equation}\label{bdSn1i}
\begin{aligned}
\|S_{n,1,i} C_{n+1} \|^2_{L^2(\DD^n)}
&\le \frac{C_T n}{N}\int_{\DD^{n+1}}
(1+\frac{|v_{i,*}|^2}{|x_{i,*}|^\theta}) |C_{n+1}|^2
\; d\mu_{[n]}d\mu_*
\\
&\le \frac{C_T}{nN} \int_{\DD^{n+1}}
\sum_{k\not =\ell}^{n+1}(1+\frac{|v_{k,\ell}|^2}{|x_{k,\ell}|^\theta}) |C_{n+1}|^2
\; d\mu_{[n]}d\mu_*
\\&\le \frac{C_T}{nN}  \| \sqrt{w_{n+1}} C_{n+1}\|_{L^2(\DD^{n+1})}^2.
\end{aligned}
\end{equation}
Similarly, we bound 
$$
\begin{aligned}
\|S_{n,0,i} C_{n}\|_{L^2(\DD^n)}^2 & \le 
\frac{C_T}{N^2}\sum_{j = 1, \; j \not =i}^n \int_{\DD^{n+1}}
(1+\frac{|v_{i,*}|^2}{|x_{i,*}|^\theta} ) |C_n(z_{[n]\setminus\{j\}}, z_*)|^2
\; d\mu_{[n]}d\mu_*
\\&\le \frac{C_T}{nN^2}  \| \sqrt{w_{n}} C_{n}\|_{L^2(\DD^{n})}^2.
\end{aligned}$$

Now recalling $D_{n,0} = \sum_i^n \nabla_{v_i}\cdot S_{n,0,i}$ and using the integration by parts similarly to \eqref{eq:bddvS-0} we have 
$$
\langle D_{n,0}C_n,C_n\rangle_{L^2_z}
=-\sum_i\langle S_{n,i,0}C_n,\nabla_{v_i}C_n\rangle_{L^2_z}
-\sum_i\langle S_{n,i,0}C_n,C_n\nabla_{v_i}\log f(z_i)\rangle_{L^2_z}.
$$
Using the Cauchy--Schwarz inequality, and also the assumption $f^{\epsilon} |\nabla_v\log f|\lesssim 1$ and H\"older's inequality as in \eqref{eq:Holder-logv}, we find that 
\begin{align*}
|\langle D_{n,0}C_{n}, C_n\rangle_{L^2_z}| &\le   \Big( \sum_i^n \|S_{n,0,i} C_{n}\|_{L^2(\DD^n)}^2 \Big)^{1/2} \Big[ \Big( \sum_i^n \|\nabla_{v_i}C_n\|_{L^2}^2\Big)^{1/2}  +   \Big( \sum_i^n \| C_n\nabla_{v_i}\log f(z_i) \|_{L^2}^2\Big)^{1/2} \Big] \\
&\le  \frac{C_T}{N}  \| \sqrt{w_{n}} C_{n}\|_{L^2(\DD^{n})} \Big[ \Big( \sum_i^n \|\nabla_{v_i}C_n\|_{L^2}^2\Big)^{1/2}  + \sqrt{n}(N/n)^{\epsilon}\Big]. 
\end{align*}
Similarly, we have
\begin{align*}
|\langle D_{n,1}C_{n+1}, C_n\rangle_{L^2_z}| \le  \frac{C_T}{N^{1/2}}  \| \sqrt{w_{n+1}} C_{n+1}\|_{L^2(\DD^{n})}  \Big[  \Big( \sum_i^n \|\nabla_{v_i}C_n\|_{L^2}^2\Big)^{1/2}+  \sqrt{n}(N/n)^{\epsilon} \Big]. 
\end{align*}
Alternatively, if we use \eqref{bdSprior}, we then obtain 
$$
\begin{aligned}
|\langle D_{n,0}C_{n}, C_n\rangle_{L^2_z}| &\le C_T n^{1/2}(N/n)^{-1+1/d} \Big[  \Big( \sum_i^n \|\nabla_{v_i}C_n\|_{L^2}^2\Big)^{1/2}+  \sqrt{n}(N/n)^{\epsilon} \Big], 
\\
|\langle D_{n,1}C_{n+1}, C_n\rangle_{L^2_z}|  &\le C_T n^{1/2}(N/n)^{-1/2+1/d} \Big[  \Big( \sum_i^n \|\nabla_{v_i}C_n\|_{L^2}^2\Big)^{1/2}+  \sqrt{n}(N/n)^{\epsilon} \Big]. 
\end{aligned}$$
The advantage of these latter estimates is that the bounds use the a priori estimates, for which we will not need to bootstrap. Finally, recalling \eqref{defDnstar}, we compute 
$$ 
\begin{aligned}
|\langle D_{n,-}C_{n-1}, C_n\rangle_{L^2_z}| \le \frac{2}{\sqrt{nN}} \sum_{i\not = j}^n \int_{\DD^n} |(K[f](x_i) - K_{i,j})| |\nabla_{v_i}C_{n-1}(z_{[n]\setminus\{j\}})| |C_n|\; d\mu_{[n]}.  
\end{aligned}
$$
Since $\nabla_{v_i}C_{n-1}$ does not depend on $z_j$, we first compute the integration $d\mu_j$, following the calculation as done above, see \eqref{bdKoneint}-\eqref{bdSn1i}. This yields 
$$ 
\Big\|  \int_{\DD} |(K[f](x_i) - K_{i,j})| |C_n|\; d\mu_j\Big\|_{L^2(\DD^{n-1})} 
\le C_Tn^{-1} \| \sqrt{w_n}C_n\|_{L^2(\DD^n)}.$$
Therefore, we bound 
$$ 
\begin{aligned}
|\langle D_{n,-}C_{n-1}, C_n\rangle_{L^2_z}| 
& \le \frac{C_T}{\sqrt{n^3N}} \sum_{i\not = j}^n \|\nabla_{v_i}C_{n-1}\|_{L^2(\DD^{n-1})}  \| \sqrt{w_n}C_n\|_{L^2(\DD^n)}
\\& \le \frac{C_T}{\sqrt{N}} \Big( \sum_{i}^n \|\nabla_{v_i}C_{n-1}\|^2_{L^2(\DD^{n-1})} \Big)^{1/2} \| \sqrt{w_n}C_n\|_{L^2(\DD^n)}
.\end{aligned}
$$
Alternatively, we may also bound 
$$ 
\Big\|  \int_{\DD} |(K[f](x_i) - K_{i,j})| |C_n|\; d\mu_j\Big\|_{L^2(\DD^{n-1})} 
\le C_T (N/n)^{1/d},$$
and therefore
$$ 
\begin{aligned}
|\langle D_{n,-}C_{n-1}, C_n\rangle_{L^2_z}|
& \le C_T (N/n)^{-1/2+1/d}\Big( \sum_{i}^n \|\nabla_{v_i}C_{n-1}\|^2_{L^2(\DD^{n-1})} \Big)^{1/2} 
.\end{aligned}
$$
The lemma thus follows. 
\end{proof}

\subsection{Order-one correlations}
In this section, we treat the order-one terms that appear in the cumulant dynamics \eqref{fullCn}-\eqref{Rnstructure}, namely those that are not treated via the a priori estimates from Lemma \ref{lem-apriori}. However, we obtain the following. 

\begin{lemma}\label{lem-oneC0}
Let $T_{n,0,1}$ and $T_{n,2}$ be defined as in \eqref{eqsCNn} and \eqref{orderone}. Then, there hold 
\begin{equation}\label{L2Tn2C}
\begin{aligned}
\|T_{n,0,1}C_n \|_{L^2(\DD^n)} &\le C_T n\|C_n\|_{L^2(\DD^{n})}, 
\\
\|T_{n,2}C_{n+2} \|_{L^2(\DD^n)} &\le C_T \| \sqrt{w_{n+2}}C_{n+2}\|_{L^2(\DD^{n+2})}.
\end{aligned}
\end{equation}
\end{lemma}
\begin{proof} Note that $a_{n,2}\le \sqrt{(n+1)(n+2)}\le 2n$. Recall that $V_{*,*'} = (K_{*,*'} - K[f](x_*))\cdot \nabla_{v_*} \log f(z_*)$, with $K_{*,*'} \sim |x_{*,*'}|^{-\alpha -1}$. Hence, 
under the assumption $\alpha< (d-1)/2$ and $\max\{0, 2\alpha+2-d\}<\theta<1$, we have \eqref{intro-bdVstar}. 
Consequently, we obtain \eqref{eq:Tn+2-intro}. 

On the other hand, we bound 
$$ \begin{aligned}
| T_{n,0,1}C_n| 
& \le \sum_{j=1}^n \int_{\DD} |V_{*,j}C_n(z_{[n]\setminus\{j\}}, z_*)|\; d\mu_* 
\\
& \le  \Big(\sum_{j=1}^n \int_{\DD} |x_{*,j}|^a|V_{*,j}|^2 \; d\mu_* \Big)^{1/2}\Big(\sum_{j=1}^n  \int_{\DD} |x_{*,j}|^{-a}|C_n(z_{[n]\setminus\{j\}}, z_*)|^2 \; d\mu_* \Big)^{1/2}
\\
&\le C_T 
\sqrt n \Big(\sum_{j=1}^n\int_{\DD} |x_{*,j}|^{-a}|C_n(z_{[n]\setminus\{j\}}, z_*)|^2 \; d\mu_* \Big)^{1/2},\end{aligned}
 $$
provided that $a>2(\alpha+1)-d$ so that $\int_{\DD} |x_{*,j}|^a|V_{*,j}|^2 \; d\mu_*\le C_T$. Therefore, we bound 
$$
\begin{aligned}
\| T_{n,0,1}C_n\|^2_{L^2}
&\le C_T n \sum_{j=1}^n \int_{\DD^{n+1}} |x_{*,j}|^{-a}|C_n(z_{[n]\setminus\{j\}}, z_*)|^2 \; d\mu_{[n]}d\mu_* 
\\
&\le C_T n \sum_{j=1}^n \int_{\DD^{n}} |C_n(z_{[n]\setminus\{j\}}, z_*)|^2 \; d\mu_{[n]\setminus\{j\}}d\mu_*
\\&\le C_T n^2 \| C_n\|_{L^2(\DD^n)}^2,
 \end{aligned}
$$
in which we have used $\sup_{z_*}\int_{\DD_j} |x_{*,j}|^{-a} \; d\mu_j\le C_T$, provided that $a<d$ (and $a>2(\alpha+1)-d$). Note that such an $a$ exists, since $\alpha < \frac{d-1}{2}$.  This completes the proof of the lemma. 
\end{proof}

\begin{proof}[Conclusion of Proposition \ref{prop-main}]
Collecting the estimates in the previous sections yields the bound on $\langle (\partial_t + L_n)C_n, C_n \rangle_{L^2_z}$. The claimed estimate for $\langle \nabla_{v_i}\cdot S_i,\varphi_n C_n\rangle_{L^2_z}$ can be obtained similarly, using the integration by parts 
$$
\left\langle \nabla_{v_i}\cdot S_i,\varphi_n C_n\right\rangle
=-\langle S_i,\varphi_n\nabla_{v_i}C_n\rangle\\
-\langle S_i,C_n\nabla_{v_i}\varphi_n\rangle\\
-\langle S_i,\varphi_nC_n\nabla_{v_i}\log f\rangle
$$
and the fact that 
$$
\left(\sum_i\|C_n\nabla_{v_i}\varphi_n\|_{L^2_z}^2\right)^{1/2}
\le C_Tn^{3/2}
$$
since $\|\varphi_n\|_\infty\le C_Tn^2$ and $\|\nabla_{v_i}\varphi_n\|_\infty\le C_Tn$ on the support of $d\mu$. 
\end{proof}

\section{Closing of cumulant bounds} \label{sec:final-cumulant}

In this section, we deduce uniform bounds on cumulants. In particular, we introduce an iterative scheme to provide close estimates on cumulants, and obtain  the following proposition. 

\begin{proposition}\label{prop-growth} Let $d\ge 2$. 
Let $C_n=C_{N,n}$ be the adjoint cumulants as constructed in Proposition \ref{prop-CNn}, satisfying the a priori bounds in \eqref{L2bounds-apriori} with a universal constant $C_0$. Assume that the source terms $g_n$ in  \eqref{eqsCNn} satisfy $\sum_{n\ge 1}\|g_n\|_{L^2} \le C_g <\infty$. In addition, assume that $\nu_N = N^{-\delta}$ for some constant $\delta \in (1/2,1)$. Then, for all $n\ge 1$, there hold 
\begin{equation}\label{key}
\sup_{t\in [0,T]}\|C_n(t)\|_{L^{2}}^2 +\nu_N\sum_i^n\|\nabla_{v_i}C_n\|^2_{L^{2}_{t,z}} + \| \sqrt{w_n} C_n\|_{L^2_{t,z}}^2 \le C_{\alpha,T} (1+C_g) n^{2a}
\end{equation}
for all $n\ge 1$, provided that $a\ge 4$ and $(1-\delta) (1/2+a/4-1) >1/d$.  
\end{proposition}

We stress that Proposition \ref{prop-growth} does not provide decay, but growth in $n$. However, the main point is of course the $N$-independence of the upper bound.  
Comparing to the a priori bounds from Lemma \ref{lem-apriori}, the above provides bounds on the weighted $L^2$ norms, which are highly non trivial. 

\begin{proof}[Proof of Proposition \ref{prop-growth}] The proposition is a combination of the previous linear analysis and iterative bounds established in the previous sections. Indeed, we start from the energy estimates in Lemmas \ref{lem-L2basic} and \ref{lem-keyL2}, namely 
\begin{equation}\label{reL2basic}
\begin{aligned}
\sup_{t\in [0,T]}\|C_n(t)\|_{L^{2}}^2 +\nu_N\sum_i^n\|\nabla_{v_i}C_n\|^2_{L^{2}_{t,z}}
&\le 2\int_0^T |\langle R_n, C_n\rangle_{L^{2}_{z}}| dt   + C_T\nu_N n 
\\&\quad + C_T (N/n)^{-\delta_\alpha } \|\sqrt{w_n}C_n\|_{L^2_{t,z}},
\end{aligned}
\end{equation}
and
 \begin{equation}\label{rekeyL2}
\begin{aligned}
c_\theta \| \sqrt{w_n} C_n\|_{L^2_{t,z}}^2 
&\le 
-2\langle R_n, \varphi_n C_n\rangle_{L^2_{t,z}} +
C_Tn^2(N/n)^{-\delta_\alpha } \|\sqrt{w_n}C_n\|_{L^2_{t,z}} +
C_T\nu_Nn^3(N/n)^{2\varepsilon}.
\\&\quad +C_Tn^2\Big[ \|C_n(T)\|_{L^2}^2+(1+\nu_N n)\|C_n\|_{L^{2}_{t,z}}^2\Big]+ \nu_Nn^2 \sum_i^n \|\nabla_{v_i}C_n\|_{L^{2}_{t,z}}^2,
\end{aligned}
\end{equation}
in which $R_n = (\partial_t + L_n)C_n$. We now apply Proposition \ref{prop-main} to bound $R_n$ that appears on the right-hand side of \eqref{reL2basic}-\eqref{rekeyL2}. Namely, we have
\begin{align}\label{rebdRn}
|\langle R_n, C_n \rangle_{L^2_z}|& \le C_T\Big[\| g_n\|_{L^2}+n (N/n)^{-1/2 + 1/d} + \nu_N \sqrt{Nn}\Big] \| C_n\|_{L^2}
\nonumber\\&\quad + \frac{C_T}{\sqrt{N}} \Big( \| \sqrt{w_{n}} C_{n}\|_{L^2(\DD^{n})} + \| \sqrt{w_{n+1}} C_{n+1}\|_{L^2(\DD^{n})} \Big) \times \nonumber\\
&\qquad \times \Big[ \Big( \sum_i^n \|\nabla_{v_i}C_n\|_{L^2}^2\Big)^{1/2} + \Big( \sum_{i}^n \|\nabla_{v_i}C_{n-1}\|^2_{L^2(\DD^{n-1})} \Big)^{1/2} +  \sqrt{n}(N/n)^{\epsilon}\Big]
\nonumber\\&\quad + C_T \Big(n\|C_n\|_{L^2(\DD^{n})}
+ \| \sqrt{w_{n+2}}C_{n+2}\|_{L^2(\DD^{n+2})}\Big) \|C_n\|_{L^2},
\end{align}
and a similar bound holds for $\langle R_n, \varphi_nC_n \rangle_{L^2_z}$ upon replacing $C_T$ by $C_Tn^2$. Alternatively, using \eqref{bdCnminus1}, we also have
\begin{align}\label{rebdRn1}
|\langle R_n, C_n \rangle_{L^2_z}|& \le C_T\Big[\| g_n\|_{L^2}+n (N/n)^{-1/2 + 1/d} + \nu_N \sqrt{Nn}\Big] \| C_n\|_{L^2}
\nonumber\\&\quad + C_T n^{1/2}(N/n)^{-1/2+1/d} \times \nonumber \\
&\qquad\times\Big[ \Big( \sum_i^n \|\nabla_{v_i}C_n\|_{L^2}^2\Big)^{1/2} +\Big( \sum_{i}^n \|\nabla_{v_i}C_{n-1}\|^2_{L^2(\DD^{n-1})} \Big)^{1/2} + \sqrt{n}(N/n)^{\epsilon} \Big]
\nonumber\\&\quad + C_T \Big(n\|C_n\|_{L^2(\DD^{n})}
+ \| \sqrt{w_{n+2}}C_{n+2}\|_{L^2(\DD^{n+2})}\Big) \|C_n\|_{L^2},
\end{align}
in which only terms in the second and third lines on the right were modified, comparing to those in \eqref{rebdRn}. We shall use the bounds \eqref{rebdRn} if $n\le (\nu_NN)^{1/4}$ and use the bounds \eqref{rebdRn1} if $n\ge (\nu_NN)^{1/4}$.  

To proceed, fix a large constant $a$, and introduce the following bootstrap quantity 
\begin{equation}\label{defbigA}
\cE_{N,T} = \sum_{n\ge 1} n^{-2a} \zeta_{n,T}^2
\end{equation}
with 
\begin{equation}\label{defzetanT}
\zeta_{n,T} = \Big[2 n^2 \sup_{t\in [0,T]}\|C_n(t)\|_{L^{2}}^2 +2\nu_Nn^2\sum_i^n\|\nabla_{v_i}C_n\|^2_{L^{2}_{t,z}} + c_\theta \| \sqrt{w_n} C_n\|_{L^2_{t,z}}^2 \Big]^{1/2}.
\end{equation}
Observe that both $\zeta_{n,T}$ and the summation $\cE_{N,T}$ are well-defined and finite, for finite $N\ge 1$, although they may blow up as $N\to \infty$. We shall prove that both $\zeta_{n,T}$ and $\cE_{N,T}$ remain finite uniformly in $N$. 

For each fixed $n$, we add up \eqref{reL2basic} and \eqref{rekeyL2}, and use the a priori bounds $\|C_n\|_{L^2}\lesssim 1$ from Lemma \ref{lem-apriori}, giving 
$$
\begin{aligned}
\zeta_{n,T}^2 &\le 2n^2\int_0^T |\langle R_n, C_n\rangle_{L^{2}_{z}}| dt +2\int_0^T |\langle R_n, \varphi_n C_n\rangle_{L^2_{z}} | dt 
\\&\quad + C_Tn^2(N/n)^{-\delta_\alpha } \|\sqrt{w_n}C_n\|_{L^2_{t,z}}
+C_Tn^3+ \nu_Nn^2 \sum_i^n \|\nabla_{v_i}C_n\|_{L^{2}_{t,z}}^2.
\end{aligned}
$$
Here we have used $\nu_N=N^{-\delta}$ and chosen $2\varepsilon\le 1/2<\delta$ to bound $\nu_N n^3(N/n)^{2\varepsilon}\le n^3$. In the above estimate, we may bound $\nu_Nn^2 \sum_i^n \|\nabla_{v_i}C_n\|_{L^{2}_{t,z}}^2 \le \frac12\zeta_{n,T}^2$, which can be absorbed into the left-hand side. In addition, we also bound $\|\sqrt{w_n}C_n\|_{L^2_{t,z}} \le c_\theta^{-1}\zeta_{n,T}$. This yields 
$$
\begin{aligned}
\zeta_{n,T}^2 &\le 4n^2\int_0^T |\langle R_n, C_n\rangle_{L^{2}_{z}}| dt  +4\int_0^T |\langle R_n, \varphi_n C_n\rangle_{L^2_{z}} | dt  + 2C_{\theta,T}n^2(N/n)^{-\delta_\alpha } \zeta_{n,T} + C_Tn^3.
\end{aligned}
$$
Here and in what follows, $C_{\theta,T}$ denotes constants that depend on $C_T$ and $c_\theta$. Next, multiplying the above inequality by $n^{-2a}$ and summing up the results to deduce bounds on $\cE_{N,T}$, see \eqref{defbigA}, we obtain  
\begin{equation}\label{bigsum}
\begin{aligned}
\cE_{N,T} &\le  \sum_{n\ge 1} n^{-2a} \Big[ 4n^2 \int_0^T |\langle R_n, C_n\rangle_{L^{2}_{z}}| dt  +4\int_0^T |\langle R_n, \varphi_n C_n\rangle_{L^{2}_{z}}| dt \Big]
\\& \quad + 
C_T\sum_{n\ge 1}n^{3-2a} + 2C_{\theta,T} \sum_{n\ge 1}n^{2-2a}(N/n)^{-\delta_\alpha } \zeta_{n,T} \end{aligned}
\end{equation}
The second summation is finite, since $a\ge 4$. On the other hand, using $\sum_n a_nb_n \le ( \sum_n a_n^2)^{1/2} (\sum_n b_n)^{1/2}$, we bound 
$$ 
\begin{aligned}
\sum_{n\ge 1} n^{2-2a}(N/n)^{-\delta_\alpha } \zeta_{n,T} \le N^{-\delta_\alpha} \Big(\sum_{n\ge 1} n^{4+2\delta_\alpha - 2a} \Big)^{1/2} \cE_{N,T}^{1/2} \le C_0N^{-\delta_\alpha} \cE_{N,T}^{1/2},
\end{aligned}
$$
which is again finite, since $\delta_\alpha <1$ and $a\ge 4$. It remains to study the first summation in \eqref{bigsum}. 

We split the summation into two cases: $n\le (\nu_NN)^{1/4}$ and $n\ge (\nu_NN)^{1/4}$. In the former case, we use \eqref{rebdRn} to deduce  
$$
\begin{aligned}
&4n^2|\langle R_n, C_n\rangle_{L^{2}_{t,z}}| +4|\langle R_n, \varphi_n C_n\rangle_{L^2_{t,z}}| 
\\& \le C_Tn^2\Big[\| g_n\|_{L^2}+n (N/n)^{-1/2 + 1/d} + \nu_N \sqrt{Nn}\Big]+C_{\theta,T} n\nu_N^{-1/2}N^{-1/2} \zeta_{n,T}\zeta_{n-1,T}
\\&\quad + C_{\theta,T} n\nu_N^{-1/2}\Big( N^{-1} \zeta_{n,T} +N^{-1/2} \zeta_{n+1,T}\Big) \zeta_{n,T}
+ C_T n^3
+ C_Tn^2\zeta_{n+2,T}\\
&\quad+C_T n^{5/2} N^{-1/2}  (N/n)^{\epsilon} (\zeta_{n,T} + \zeta_{n+1,T}),
\end{aligned}
$$
in which $C_{\theta,T}$ denotes constants that depend on $C_T$ and $c_\theta$. Using $n\le (\nu_NN)^{1/4}$, we may further simplify the above bounds to get 
\begin{align}\label{bd4Rn1}
&4n^2|\langle R_n, C_n\rangle_{L^{2}_{z}}| +4|\langle R_n, \varphi_n C_n\rangle_{L^2_{z}}| \nonumber
\\& \le C_Tn^2\Big[\| g_n\|_{L^2}+n (N/n)^{-1/2 + 1/d} + \nu_N \sqrt{Nn}\Big]+C_{\theta,T} (\nu_NN)^{-1/4} \zeta_{n,T}\zeta_{n-1,T}
\\&\quad + C_{\theta,T} (\nu_NN)^{-1/4} \Big(\zeta_{n,T} +\zeta_{n+1,T}\Big) \zeta_{n,T}
+ C_T n^3
+ C_Tn^2\zeta_{n+2,T} + C_T n^{5/2} N^{-1/2}  (N/n)^{\epsilon} (\zeta_{n,T} + \zeta_{n+1,T}). \nonumber
\end{align}
Therefore, taking the summation in $1\le n\le (\nu_NN)^{1/4}$ and using 
\begin{align} \label{eq:add-term-ENT}
N^{-1/2}
\sum_{n\ge1}n^{5/2-2a} (N/n)^{\epsilon}
(\zeta_{n,T}+\zeta_{n+1,T})
\le C_{a,\epsilon}E_{N,T}^{1/2},
\end{align}
with $a\ge 4$, we can  bound 
$$
\begin{aligned}
&
\sum_{1\le n\le (\nu_NN)^{1/4}}n^{-2a}\Big[4n^2|\langle R_n, C_n\rangle_{L^{2}_{z}}| +4|\langle R_n, \varphi_n C_n\rangle_{L^2_{z}}| \Big]
\\& \le C_T\sum_{n\ge 1}n^{2-2a}\Big[\| g_n\|_{L^2}+n (N/n)^{-1/2 + 1/d} + \nu_N \sqrt{Nn}\Big]+C_{\theta,T} (\nu_NN)^{-1/4} \sum_{n\ge 1}n^{-2a}\zeta_{n,T}\zeta_{n-1,T}
\\&\quad + C_{\theta,T} (\nu_NN)^{-1/4} \sum_{n\ge 1}n^{-2a}\Big(\zeta_{n,T} +\zeta_{n+1,T}\Big) \zeta_{n,T}
+ C_T \sum_{n\ge 1}n^{3-2a}
+ C_T\sum_{n\ge 1}n^{2-2a}\zeta_{n+2,T} + C_{a,\epsilon}E_{N,T}^{1/2}. 
\end{aligned}
$$
We bound 
$$ \sum_{n\ge 1}n^{2-2a}\Big[\| g_n\|_{L^2}+n (N/n)^{-1/2 + 1/d} + \nu_N \sqrt{Nn}\Big]  \le C_0 \Big( 1+ \nu_N \sqrt{N}+ \sum_{n\ge 1} \| g_n\|_{L^2}\Big).$$
On the other hand, using again $\sum_n a_nb_n \le ( \sum_n a_n^2)^{1/2} (\sum_n b_n)^{1/2}$, we bound 
$$
\begin{aligned}
(\nu_NN)^{-1/4} \sum_{n\ge 1}n^{-2a}\zeta_{n,T}\zeta_{n-1,T} &\le (\nu_NN)^{-1/4} \cE_{N,T}, 
\\(\nu_NN)^{-1/4} \sum_{n\ge 1}n^{-2a}\Big(\zeta_{n,T} +\zeta_{n+1,T}\Big) \zeta_{n,T} &\le (\nu_NN)^{-1/4} \cE_{N,T},
\\
\sum_{n\ge 1}n^{2-2a}\zeta_{n+2,T} &\le C_0 \cE_{N,T}^{1/2}.
\end{aligned}
$$
Therefore, the sum over $1\le n\le (\nu_NN)^{1/4}$ is bounded by 
\begin{equation}\label{bdRfirst}
\begin{aligned}
&\sum_{1\le n\le (\nu_NN)^{1/4}}n^{-2a}\Big[4n^2 |\langle R_n,  C_n\rangle_{L^2_{z}}| +4|\langle R_n, \varphi_n C_n\rangle_{L^2_{z}}| \Big] 
\\& \le C_T (1+\nu_N\sqrt{N}+ \sum_{n\ge 1} \| g_n\|_{L^2}) + C_{\theta,T} (\nu_NN)^{-1/4} \cE_{N,T} + C_T \cE_{N,T}^{1/2}.
\end{aligned}
\end{equation}
Finally, we consider the case when $n\ge (\nu_NN)^{1/4}$. In this case, we use instead the following bounds 
\begin{align}\label{bd4Rn2}
&4n^2|\langle R_n, C_n\rangle_{L^{2}_{z}}| +4|\langle R_n, \varphi_n C_n\rangle_{L^2_{z}}| 
\nonumber\\& \le C_Tn^2\Big[\| g_n\|_{L^2}+n (N/n)^{-1/2 + 1/d} + \nu_N \sqrt{Nn}\Big]+C_{\theta,T} n^2\nu_N^{-1/2}(N/n)^{-1/2+1/d} \zeta_{n-1,T}
\\&\quad + C_{\theta,T} n^2\nu_N^{-1/2} (N/n)^{-1/2+1/d} \zeta_{n,T}
+ C_T n^3
+ C_Tn^2\zeta_{n+2,T} + C_T n^{5/2} N^{-1/2}  (N/n)^{\epsilon} (\zeta_{n,T} + \zeta_{n+1,T}).\nonumber
\end{align}
Taking the summation in $n\ge (\nu_NN)^{1/4}$, with weight $n^{-2a}$, the first and the last two terms are estimated exactly the same as before. As for the remaining terms, we bound 
$$
\begin{aligned}
\sum_{n\ge (\nu_NN)^{1/4}} n^{2-2a}\nu_N^{-1/2}(N/n)^{-1/2+1/d} \zeta_{n-1,T} 
&\le N^{-1/2+1/d} \nu_N^{-1/2}\Big( \sum_{n\ge (\nu_NN)^{1/4}} n^{4-2a}n^{1-2/d} \Big)^{1/2} \cE_{N,T}^{1/2}
\\
&\le C_a N^{1/d}(\nu_NN)^{(1-1/d-a)/4}\cE_{N,T}^{1/2}
\end{aligned}
$$
provided that $a \ge 4$. The term $C_{\theta,T} n^2\nu_N^{-1/2} (N/n)^{-1/2+1/d} \zeta_{n,T}$ is estimated similarly. 
The terms involving $\sqrt {n}(N/n)^\epsilon$ in Lemma \ref{lem-Cnminus} are estimated using \eqref{eq:add-term-ENT} where the sum on the left-hand side is taken for all $n\ge 1$. 
Therefore, combining the two cases, we have obtained 
\begin{equation}\label{bdR2nd}
\begin{aligned}
&\sum_{n\ge 1}n^{-2a}\Big[4n^2|\langle R_n, C_n\rangle_{L^{2}_{z}}| +4|\langle R_n, \varphi_n C_n\rangle_{L^2_{z}}| \Big] 
\\& \le C_T (1+\nu_N\sqrt{N}+ \sum_{n\ge 1} \| g_n\|_{L^2}) + C_{\theta,T} (\nu_NN)^{-1/4} \cE_{N,T} 
\\&\quad + C_{\theta,T} \Big(1+ N^{1/d}(\nu_NN)^{(1-1/d-a)/4}\Big)\cE_{N,T}^{1/2}.
\end{aligned}
\end{equation}
Putting this above estimate back into \eqref{bigsum}, we obtain 
\begin{equation}\label{bigcE1}
\begin{aligned}
\cE_{N,T} & \le C_T \Big(1+\nu_N\sqrt{N}+ \sum_{n\ge 1} \| g_n\|_{L^2}\Big) + C_{\theta,T} (\nu_NN)^{-1/4} \cE_{N,T} 
\\&\quad + C_{\theta,T} \Big(1+ N^{1/d}(\nu_NN)^{(1-1/d-a)/4} + N^{-\delta_\alpha}\Big)\cE_{N,T}^{1/2}.
\end{aligned}
\end{equation}
Under the assumption that $\nu_N = N^{-\delta}$ with $\delta \in (1/2,1)$, we note that $\nu_N \sqrt{N} \ll 1 $ and $\nu_N N \gg1$. This latter condition allows us to absorb the second term on the right into the left-hand side. In addition, we take $a$ large enough so that $(1-\delta) (1/2+a/4-1) >1/d$, which implies that 
$N^{1/d}(\nu_NN)^{(1-1/d-a)/4} \le 1$. This yields 
$$
\cE_{N,T}  \le C_T \Big(1+ \sum_{n\ge 1} \| g_n\|_{L^2}\Big) + o(1) \cE_{N,T} + C_T\cE_{N,T}^{1/2} 
$$
which implies that
\begin{equation}\label{bigcE2}
\begin{aligned}
\cE_{N,T} & \le C_T \Big(1+ \sum_{n\ge 1} \| g_n\|_{L^2}\Big) .
\end{aligned}
\end{equation}
The proposition thus follows, upon recalling \eqref{defbigA}-\eqref{defzetanT}.
\end{proof}
 
\section{Proof of Theorem \ref{theo-main}}\label{sec:final-uniqueness}
In this section, we conclude   the main result. From the uniform bounds established in Proposition \ref{prop-growth}, up to an extraction of subsequences of $C_{N,n}$, we obtain the weak convergence 
$$ C_{N,n} \rightharpoonup  C_{\infty,n} , \qquad \sqrt{\nu_N} \nabla_{v_i}C_{N,n} \rightharpoonup  W_{\infty,n} $$
in $L^2([0,T]\times \DD^2)$, as $N\to \infty$. In addition, $C_{N,n} \rightharpoonup  C_{\infty,n}$ converges weakly in weighted $L^2$ spaces with weight $w_n$. As a result, upon taking the limit of $\nu_N \sqrt{N} \to 0$ and $N\to \infty$ in the hierarchy \eqref{eqsCNn}, the limiting cumulants $C_{\infty,n}$ solve
\begin{equation}\label{limitCn1}
(\partial_t + \cL^{\otimes n}) C_{\infty,n} = T_{n,0}C_{\infty,n} + T_{n,2}C_{\infty,n+2}
\end{equation}   
with zero initial data $C_{\infty,n}(0)=0$, where $\cL = v\cdot \nabla_x + K[f](x)\cdot \nabla_v$ denotes the one-particle Vlasov operator \eqref{VP}, and 
\begin{equation}\label{limitTn1}
\begin{aligned}T_{n,0}C_{\infty,n} &= \sum_{j=1}^n \int_{\DD} V_{*,j}C_{\infty,n}(z_{[n]\setminus\{j\}}, z_*)\; d\mu_* ,  
\\T_{n,2}C_{\infty,n+2} &=\sqrt{(n+1)(n+2)}\iint_{\DD^2} V_{*,*'} C_{\infty,n+2}(z_{[n]}, z_*, z_*')
\; d\mu_*d\mu_{*}'.
\end{aligned}
\end{equation}
Here $V_{1,2} = (K_{1,2} - K[f](x_1))\cdot \nabla_{v_1} \log f(z_1)$. 

Let us now end the proof of Theorem \ref{theo-main}, assuming the uniqueness of \eqref{limitCn1} which will be proved below.  Indeed, we recall from Corollary \ref{cor-main} that 
\begin{equation}\label{goalagain}
\begin{aligned}
&
\langle F_{N,k} - f^{\otimes k}, \psi^{\otimes k}\rangle_{L^2([0,T]\times \DD^k)} 
\\&= N \int_0^T \int_{\DD^2} V_{1,2}\Phi_{N,2} \; d\mu_1d\mu_2  - \nu_NN \int_0^T \int_{\DD}B_1\Phi_{N,1} \; d\mu_1 -  \langle F_N(0) - f_0^{\otimes N}, \Phi_N(0) \rangle_{L^2(\DD^N)}.
\end{aligned}
\end{equation}
We may assume that $F_N(0) = f_0^{\otimes N}$, and so the last term in \eqref{goalagain} vanishes. Next, by definition, we note that 
$$ \nu_NN \int_0^T \int_{\DD}B_1\Phi_{N,1} \; d\mu_1 = \nu_N \sqrt{N} \int_0^T \int_{\DD}B_1C_{N,1} \; d\mu_1 \to 0$$ 
as $N\to \infty$, since $C_{N,1}$ is uniformly bounded in $L^2$ and $\nu_N \sqrt{N} \to 0$ by assumption. Finally, using \eqref{intro-bdVstar} we note that  
$$
\begin{aligned}
N \int_0^T \int_{\DD^2} V_{1,2}\Phi_{N,2} \; d\mu_1d\mu_2
& = \sqrt{\frac{2N}{N-1}}
\int_0^T \int_{\DD^2} V_{1,2}C_{N,2} \; d\mu_1d\mu_2
\\
&\to \sqrt{2} \int_0^T \int_{\DD^2} V_{1,2}C_{\infty,2} \; d\mu_1d\mu_2
,\end{aligned} $$
which vanishes, provided that $C_{\infty,n}=0$. This proves that $\langle F_{N,k} - f^{\otimes k}, \psi^{\otimes k}\rangle_{L^2([0,T]\times \DD^k)} \to 0$ in view of \eqref{goalagain}, and thus completes the proof of Theorem \ref{theo-main}.

Thus to conclude Theorem \ref{theo-main} it remains to prove that the limiting hierarchy \eqref{limitCn1} has a unique solution. More precisely, we need to prove that \eqref{limitCn1} only has trivial solution $C_{\infty,n}(t)=0$ provided that the initial data is zero.  For convenience, we write the hierarchy \eqref{limitCn1} as 
\begin{equation}\label{hierCn}
(\partial_t + \cL_n) C_{\infty,n} = \sqrt{(n+1)(n+2)} T_{n+2} C_{\infty, n+2},
\end{equation}
where $\cL_n = \sum_{j=1}^n \cL^{(j)}$, with 
\begin{equation}\label{cLparticlej}
	\cL^{(j)} C_n = v_j\cdot \nabla_{x_j}C_n+ K[f](x_j)\cdot \nabla_{v_j}C_n  - \int_{\DD} V_{*,j}C_n(z_{[n]\setminus\{j\}}, z_*)\; d\mu_*,
\end{equation}
and \begin{equation}\label{limitTn2}
	\begin{aligned}
		T_{n+2}C_{\infty,n+2} &=\iint_{\DD^2} V_{*,*'} C_{\infty,n+2}(z_{[n]}, z_*, z_*')
		\; d\mu_*d\mu_{*}'.
	\end{aligned}
\end{equation}

Note that $\cL_n=\cL_n(t)=\sum_{j=1}^n \cL^{(j)}(t)$ and $T_{n+2}=T_{n+2}(t)$ depend on time, which we sometimes omit in the notation. Since  $\{\cL^{(j)}(t)\}_j$ acts on different particle variables, the propagator for the dynamics $(\partial_t+\cL_n(t))F_n=0$ on $0\le s\le t\le T$ with the initial value $U_n(s,s)=1$ has the factorized form 
 \begin{equation}\label{eq:Un-factorized}
 U_n(t,s)=\prod_{j=1}^nU^{(j)}(t,s).
\end{equation}
where $U^{(j)}(t,s)$ is the corresponding propagator in the $j$-th variable. We have the following Dyson series expansion. 

\begin{proposition}\label{prop-Dyson} Let $(C_{\infty,n}(t))_n$ be the solution to the infinite hierarchy \eqref{hierCn} with zero initial data. Then, for any $k\ge 1$, there holds
	\begin{equation}\label{DysonCninf0}
		\begin{aligned}
			C_{\infty,n}(t) 
			&=b_{n,k} \int_0^t U_n(t,s)
			\prod_{j=1}^k\Big(\int_s^t\cB_{n+2j}(\tau,s) \; d\tau\Big) G_{n+2k}(s)\; ds,
		\end{aligned}
	\end{equation}
	where $0\le b_{n,k}  \lesssim_n 2^k k^{n/2+1}$, $G_{n+2k}(t) = T_{n+2k+2}(t)C_{\infty,n+2k+2}(t)$ and 
	$$
	\cB_{n+2j}(t,s) = T_{n+2j}(t) U^{(n+2j-1)}(t,s) U^{(n+2j)} (t,s). 
	$$
\end{proposition}
\begin{proof} 
	For any fixed $n\ge 1$, from  \eqref{hierCn} we have the Duhamel formula 
	$$
	C_{\infty,n}(t) = a_{n+2} \int_0^t U_n(t,s)  T_{n+2}C_{\infty,n+2}(s)\; ds,$$
	with $a_{n+2} = \sqrt{(n+1)(n+2)}$. 
	Iteratively, for any $k\ge 1$, we may write 
	$$
	\begin{aligned}
	C_{\infty,n}(t) 
	&= a_{n+2} a_{n+4}\int_0^t \int_0^{t_1} U_n(t,t_1) T_{n+2} U_{n+2} (t_1,t_2) T_{n+4}C_{\infty,n+4}(t_2)\; dt_2dt_1 = \cdots 
	\\
	&= \prod_{j=1}^{k+1} a_{n+2j} \int_0^t \int_0^{t_1}\cdots \int_0^{t_{k}} U_n(t,s) \cB_{n+2}(t_1,s) \cdots \cB_{n+2k}(t_{k},s) G_{n+2k}(s)\; dsdt_{k} \cdots dt_1
	\end{aligned}$$
	where the integration is taken over $\{ 0<s<t_k<\cdots<t_1<t\}$.	Here we used the semigroup property $U_n(t,s)=U_n(t,t_1) U_n(t_1,s)$ and the identity  
	$$\cB_{n+2j}(t,s) = T_{n+2j}(t)U^{(n+2j-1)}(t,s) U^{(n+2j)} (t,s) = U_{n+2j-1}(s,t) T_{n+2j}(t) U_{n+2j} (t,s)$$
	which follows from the factorized form \eqref{eq:Un-factorized} and the fact that $U^{(j)}(t,s)$  acts precisely on the $j$-th particle.

	Next, note that $G_{n+2k}(\cdot)$ is a symmetric function of $n+2k$ particles, and the propagators $\{\cB_{n+2j}(t,s)\}_{j}$ have the same form except that  each operator $\cB_{n+2j}(t,s)$ acts precisely on the last two particles $n+2j-1$ and $n+2j$, keeping the remaining variables intact. 
	As a result, 
	$$\cB_{n+2}(t_1,s) \cdots \cB_{n+2k}(t_{k},s) G(s) = \prod_j^{k}\cB_{n+2\sigma(j)}(t_{\sigma(j)},s) G(s)\circ Z_{\sigma}$$
	for any permutation $\sigma$ of $\{1,\cdots, k\}$ and $G\circ Z_\sigma = G(z_1,\cdots, z_n, z_{n+\sigma(1)},\cdots, z_{n+\sigma(k)})$, namely the permutation of the last $2k$ particles. Using the symmetry, we thus obtain  
	\begin{align*}
	&\int_{\{s<t_k<\cdots<t_1<t\}} \cB_{n+2}(t_1,s) \cdots \cB_{n+2k}(t_{k},s)G_{n+2k}(s) \; dt_{k} \cdots dt_1\\
	&= \frac{1}{k!} \int_{[s,t]^k} \cB_{n+2}(t_1,s) \cdots \cB_{n+2k}(t_{k},s) \; dt_{k} \cdots dt_1 G_{n+2k}(s)\\
&= \frac{1}{k!} \prod_{j=1}^k\Big(\int_s^t\cB_{n+2j}(\tau,s) \; d\tau\Big)G_{n+2k}(s).
	\end{align*}
Therefore, 
	\begin{equation}\label{DysonCninf}
		\begin{aligned}
			C_{\infty,n}(t) 
			&= \frac{1}{k!}\prod_{j=1}^{k+1} a_{n+2j} \int_0^t U_n(t,s)
			\prod_{j=1}^k\Big(\int_s^t\cB_{n+2j}(\tau,s) \; d\tau\Big) G_{n+2k}(s)\; ds.
		\end{aligned}
	\end{equation}
Finally, we note that by Stirling's formula, 
	\begin{equation}\label{coeffak}
	b_{n,k}=\frac{1}{k!}	\prod_{j=1}^{k+1} a_{n+2j} = \frac{1}{k!} \sqrt{\frac{(n+2k+2)!}{n!}} \lesssim_n 2^{k} k^{n/2+1} .
	\end{equation}
This ends the proof of the proposition. 
\end{proof}

\begin{lemma}\label{lem-expLnt} For all $0\le s\le t\le T$ and $n\ge 1$ we have
	\begin{equation} 
		\|U_n(t,s) G_n\|_{L^2(\DD^n)} \le e^{C_T n (t-s)} \|G_n\|_{L^2(\DD^n)}.
	\end{equation}
\end{lemma}

\begin{proof} Due to the factorized form \eqref{eq:Un-factorized} it suffices to consider the one-particle dynamics 
	\begin{equation}\label{VPBF} (\partial_t + \cL^{(1)}) F = 0,\end{equation}
where 
	\begin{equation}\label{cLparticle1}
\cL^{(1)} F = v_1\cdot \nabla_{x_1}F+ K[f](x_1)\cdot \nabla_{v_1}F  - B^{(1)}F
\end{equation}
and
\begin{equation}\label{eq:B1-def}
B^{(1)}F = \int_{\DD} V_{*,1}F(z_*)\; d\mu_*.
\end{equation}

Taking the $L^2$ inner product with respect to the measure $d\mu = f dz$, we compute 
$$ \frac12 \frac{d}{dt} \|F\|_{L^2(\DD)}^2 = \frac12 \int |F|^2 \Big( \partial_t + v \cdot \nabla_x + K[f]\cdot \nabla_v \Big) \log f \; d\mu +  \langle F, B^{(1)}F\rangle_{L^2(\DD)}.$$
The first integral term on the right vanishes, thanks to Vlasov's equation \eqref{VP}. As for the second term, recall that $V_{1,2} = (K_{1,2} - K[f](x_1))\cdot \nabla_{v_1} \log f(z_1)$ and $f^{\varepsilon} |\nabla_v \log f|\lesssim_T 1$. Therefore, 
\begin{equation}\label{eq:B1-bound} 
\|B^{(1)}F\|_{L^2(\DD)} \le C_T\|F\|_{L^2(\DD)}
\end{equation}
where we used $K\in L^1_{\rm loc}$ and Young's inequality, upon recalling that $d\mu$ is a compactly supported probability measure. Consequently, the inner product $\langle F, B^{(1)}F\rangle_{L^2(\DD)}$ is bounded by $C_T \|F\|_{L^2(\DD)}^2$. This proves that $\|F(t)\|_{L^2(\DD)} \le e^{C_T t} \|F(0)\|_{L^2(\DD)}
$ by Gr\"onwall's inequality, and the lemma thus follows. 
\end{proof}

\begin{lemma}\label{lem-bdCBt}
	Let $\cB(t,s) = T_{2}(t) U_2(t,s)$ be the two-particle contraction operator.  
	Then, for any constant $\max\{0,2+2\alpha-d\}<\theta<1$, and for $|t-s|$ sufficiently small, there holds 
	\begin{equation}
		\Big| \int_s^t\cB(\tau,s) G(s)\; d\tau\Big| \le C_T |t-s|^{1-\theta/2}\|G(s)\|_{L^2(\DD^2)}
	\end{equation}
	for any two-particle distribution $G$. 
\end{lemma}
\begin{proof} Recall that the propagator $U_2(t,s)=U^{(1)}(t,s)U^{(2)}(t,s)$ is generated by 
$$
\sum_{j=1}^2 
\Big( v_j\cdot\nabla_{x_j}
+K[f](t,x_j)\cdot\nabla_{v_j} \Big) - R(t),\qquad R =B^{(1)}+B^{(2)}.
$$
Since the nonlocal part $B^{(j)}$ given in \eqref{eq:B1-def} may complicate the time average analysis, let us first remove it using Lemma \ref{lem-expLnt}. To be precise, let $\widetilde{U}_2(t,s)$ denote
the two-particle propagator associated with the Vlasov generator (without $R(t)$), and let us assume for the moment that 
\begin{equation}\label{eq:Vlasov-bound}
\left|
\int_s^t T_2(\tau)\widetilde{U}_2(\tau,s)G\,d\tau
\right|
\le C_T(t-s)^{1-\theta/2}
\|G\|_{L^2(\DD^2)}
\end{equation}
for some $0<\theta<1$
which will be proved later. By Duhamel's formula,
\[
U_2(\tau,s)=\widetilde{U}_2(\tau,s)
+\int_s^\tau \widetilde{U}_2(\tau,\xi)
R(\xi)U_2(\xi,s)\,d\xi.
\]
we can write  
\begin{align*}
&\int_s^t T_2(\tau)U_2(\tau,s)G \,d\tau=\int_s^t T_2(\tau)\widetilde{U}_2(\tau,s)G\,d\tau+
\int_s^t d\tau \int_s^\tau d\xi 
T_2(\tau)\widetilde{U}_2(\tau,\xi)
 R(\xi)U_2(\xi,s)G. 
\end{align*}
The first term is already controlled by \eqref{eq:Vlasov-bound}. For the second term, using \eqref{eq:Vlasov-bound} with $G$ replaced by $R(\xi)U_2(\xi,s)G$, the $L^2$-boundedness of $R(\xi)$ from \eqref{eq:B1-bound}, and the bound of $U_2(\xi,s)$ in Lemma \ref{lem-expLnt}, we have 
\begin{align*}
&\left| \int_s^t d\tau \int_s^\tau d\xi 
T_2(\tau)\widetilde{U}_2(\tau,\xi)
 R(\xi)U_2(\xi,s)G\right| =\left| \int_s^t d \xi \int_\xi^t d\tau  T_2(\tau)\widetilde{U}_2(\tau,\xi)
 R(\xi)U_2(\xi,s)G  \right|  \\
 &\le \int_s^t d \xi  \left| \left( \int_\xi^t d\tau  T_2(\tau)\widetilde{U}_2(\tau,\xi) \right)
 R(\xi)U_2(\xi,s)G  \right| \lesssim_T  \int_s^t d \xi |t-\xi|^{1-\theta/2} \|  R(\xi)U_2(\xi,s)G \|_{L^2(\DD^2)}\\
 &\lesssim_T \int_s^t d \xi |t-\xi|^{1-\theta/2}  \|  G \|_{L^2(\DD^2)}\lesssim_T  |t-s|^{1-\theta/2} \|  G \|_{L^2(\DD^2)}. 
  \end{align*}
Thus, it remains only to prove
\eqref{eq:Vlasov-bound}.

	Introduce the Vlasov one-particle  characteristic $Z_{s,t} = (X_{s,t}, V_{s,t})$ where 
	$$ \dot X_{s,t} = V_{s,t}, \qquad \dot V_{s,t} = K[f](s,X_{s,t})$$
	with initial data $z=(x,v)$ at $s=t$. 
Since $d\mu = f dz$, the left-hand side of \eqref{eq:Vlasov-bound} can be reformulated as 
	\begin{align}\label{eq:Vlasov-bound-Z}
I_{s,t}&=\left| \int_s^t d\tau  \iint_{\DD^2} V_{12}(\tau,z_1,z_2)
	f(\tau,z_1)f(\tau,z_2) G(Z_{s,\tau}^{(1)}, Z_{s,\tau}^{(2)}) \;dz_1 dz_2 \right|.
	\end{align}
Recall that $V_{1,2}(\tau,z_1,z_2) = (K(x_1-x_2) - K[f](\tau,x_1))\cdot \nabla_{v_1} \log f(\tau,z_1)$ with $K[f] \in L^\infty$ and $|\nabla_v \log f| f^{\epsilon} \lesssim_{T,\epsilon} 1$ for $\epsilon>0$ small by \eqref{defmaxPT}. Using also that $z\mapsto Z_{s,t}(z)$ is a volume-preserving map and $f(\tau,Z_{\tau,s}(z))=f(s,z)$, we can bound 
\begin{align*}
I_{s,t} 
		&\lesssim  \int_s^t d\tau  \iint_{\DD^2} dz_1 dz_2 (1+ |K(x_1-x_2)|)   f(\tau,z_1)^{1-\epsilon } f(\tau,z_2) |G(Z_{s,\tau}^{(1)}, Z_{s,\tau}^{(2)})| \\
		&= \iint_{\DD^2} dz_1 dz_2 \left( \int_s^t d\tau  (1+ |K(X^{(1)}_{\tau,s}-X^{(2)}_{\tau,s})| ) \right)   f(s,z_1)^{1-\epsilon } f(s,z_2) |G(z_1,z_2)|\\
		&\le \sqrt{J_{s,t}} \|G\|_{L^2(\DD^2)}
	\end{align*}
	with 
	\begin{align}\label{bdIKst}
			J_{s,t}&= \iint_{\DD^2} \left |\int_s^t d \tau (1+ |K(X^{(1)}_{\tau,s}-X^{(2)}_{\tau,s})|) \right|^2 f(s,z_1)^{1-2\epsilon } f(s,z_2) \; dz_1 dz_2\\
			&=  \int_s^t\int_s^t\iint_{\DD^2} (1+ |K(X^{(1)}_{\tau_1,s}-X^{(2)}_{\tau_1,s})|) (1+ |K(X^{(1)}_{\tau_2,s}-X^{(2)}_{\tau_2,s})|)  f(s,z_1)^{1-2\epsilon } f(s,z_2) 
			\; dz_1 dz_2 d\tau_1d\tau_2.\nonumber
	\end{align}
	To exploit phase mixing, we write 
	\begin{equation}\label{straightX} X_{s,t} = x - (t-s)\Psi_{s,t}, \qquad \Psi_{s,t} = v -  \frac{1}{t-s}\int_s^t (\tau-s)K[f](\tau,X_{\tau,t})\; d\tau.
	\end{equation}
	Since $K[f]$ is sufficiently smooth, we have $\| \Psi_{s,t} - v\|_{L^\infty_{x,v}} \lesssim |t-s|$. In addition, we compute 
	$$
	\begin{aligned}
	\partial_x( \Psi_{s,t} -v)
	&= - \frac{1}{t-s}\int_s^t (\tau-s)\partial_x K[f](\tau,X_{\tau,t}) \Big[ 1 - (t-\tau)\partial_x\Psi_{\tau,t}\Big]\; d\tau 
	\\
	\partial_v ( \Psi_{s,t} -v)
	&=  \frac{1}{t-s}\int_s^t (\tau-s)(t-\tau)\partial_x K[f](\tau,X_{\tau,t}) \partial_v \Psi_{\tau,t}\; d\tau 
	. \end{aligned}$$
	To bound the characteristic, we let 
	$$\zeta(s,t) = \sup_{s\le \tau\le t} \Big[ \|  \partial_x ( \Psi_{\tau,t} -v)\|_{L^\infty_{x,v}} + \|  \partial_v ( \Psi_{\tau,t} -v)\|_{L^\infty_{x,v}}\Big],$$ 
	and iteratively compute 
	$$
	\begin{aligned} 
	\zeta(s,t) &\le \frac{1}{t-s}\int_s^t (\tau-s)( 1 + |t-\tau|) \|\partial_x K[f](\tau)\|_{L^\infty_x} \Big( 1 + \zeta(\tau,t)\Big)\; d\tau
	\\&\lesssim |t-s|\sup_{0\le s\le \tau\le t}\|\partial_x K[f](\tau)\|_{L^\infty_x}  (1+\zeta(s,t)). 
	\end{aligned}$$
	Therefore, for sufficiently small $|t-s|\ll1$, we may absorb $\zeta(s,t)$ into the left-hand side, yielding 
	$$
	\begin{aligned} 
	\zeta(s,t) &\lesssim |t-s|\sup_{0\le s\le \tau\le t}\|\partial_x K[f](\tau)\|_{L^\infty_x} . 
	\end{aligned}$$
	In particular, there hold 
	$$|\det \nabla_xX_{s,t}| \gtrsim 1, \qquad |\det \nabla_v \Psi_{s,t}|\gtrsim 1,$$ 
	for $|t-s|\ll1$. Thus the map $v\mapsto \Psi_{s,t}(x,v)$ is a $C^1$ diffeomorphism in $\RR^d$ and the characteristic $X_{s,t}$ can be straightened as in \eqref{straightX}. Returning to \eqref{bdIKst}, we introduce a change of variables 
	$$ y_1 = X^{(1)}_{\tau_1,s}-X^{(2)}_{\tau_1,s}$$
	with Jacobian determinant $\det(\nabla_{x_1} y_1) \gtrsim 1$. Note that 
	$$ X^{(1)}_{\tau_1,s}-X^{(1)}_{\tau_2,s} = (\tau_1-\tau_2) W^{(1)}(x_1,v_1), $$
	with 
	$$W^{(1)} (x_1,v_1)= v_1 + \frac{1}{\tau_1-\tau_2}\int_{\tau_1}^{\tau_2} (\theta - \tau_1)K[f](\theta, X^{(1)}_{\theta,s})\; d\theta - \int_{\tau_2}^s K[f](\theta, X^{(1)}_{\theta,s})\; d\theta.$$
	Proceeding similarly as the above analysis, we see that the map $(x_1,v_1) \mapsto (y_1,w_1)$ with $w_1 = W^{(1)}(x_1,v_1)$ is a $C^1$ diffeomorphism, with a 
Jacobian determinant equal to $1+O(t-s)$. 
	Therefore, by changing the variables and using the uniform bound of the corresponding Jacobian, we bound 
	$$
	\begin{aligned}
	J_{s,t}
	\lesssim \int_s^t\int_s^t\iint_{\DD^2} (1+|K(y_1)|) (1+ |K( y_1  - (\tau_1-\tau_2) (w_1 - w_2)|)  f(s,z_1)^{1-2\epsilon } f(s,z_2)  dz_1 dz_2 d\tau_1d\tau_2.
	\end{aligned}
	$$
	Now recalling $|K(x)|\lesssim  |x|^{-\alpha-1}$, we have the elementary bound 
	$$ \int_{\{|y_1|\lesssim 1\}} (1+|K(y_1)|)(1+ |K( y_1  - (\tau_1-\tau_2) (w_1 - w_2)|) \; dy_1 \lesssim 1+ |(\tau_1-\tau_2) (w_1 - w_2)|^{-\theta}$$
	for any fixed constant $\max\{0,2\alpha+2-d\}<\theta<1$. Note that since $\alpha<(d-1)/2$, the choice of $\theta$ is possible. Since the singularity in $w_1$ is integrable, we can bound 
	$$\int_s^t\int_s^t |\tau_1-\tau_2|^{-\theta} \; d\tau_1d\tau_2 \lesssim |t-s|^{2-\theta}.$$
	This proves \eqref{eq:Vlasov-bound}, and the lemma.   
\end{proof}

\begin{proof}[Proof of uniqueness]
We are now ready to prove that solutions to the limiting  hierarchy \eqref{hierCn} are unique. Precisely, for zero initial data and for any $k\ge 1$, solutions $C_{\infty,n}(t) $ to the hierarchy satisfy the Dyson series expansion \eqref{DysonCninf0}. Using Lemmas \ref{lem-expLnt} and \ref{lem-bdCBt}, we bound 
$$
		\begin{aligned}
			\|C_{\infty,n}(t)\|_{L^2(\DD^n)} 
			&\lesssim k^{n/2} 2^k \int_0^t e^{C_Tn (t-s)}
			\Big\| \prod_{j=1}^k\Big(\int_s^t\cB_{n+2j}(\tau,s) \; d\tau\Big) G_{n+2k}(s)\Big\|_{L^2(\DD^n)}\; ds
			\\
			&\lesssim  k^{n/2} 2^k  t^{k(1-\theta/2)}  \int_0^t e^{C_Tn (t-s)} \|G_{n+2k}(s)\Big\|_{L^2(\DD^{n+2k})}\; ds
		\end{aligned}
$$
for some constant $\max\{0,2+2\alpha-d\}<\theta<1$. By the Cauchy--Schwarz inequality, we have
$$
\begin{aligned}
\|G_{n+2k}(s)\|_{L^2(\DD^{n+2k})} 
 &\le \| T_{n+2k+2}C_{\infty,n+2k+2}(s)\|_{L^2(\DD^{n+2k})} 
 \lesssim \|\sqrt{w_{n+2k+2}}C_{\infty,n+2k+2}(s)\|_{L^2(\DD^{n+2k+2})}
\end{aligned}		
$$
in which the last estimate we used \eqref{L2Tn2C}, upon recalling \eqref{limitTn2}. Finally, using the a priori bounds \eqref{key}, we thus obtain 
\begin{align} \label{eq:final-Dyson}
			\|C_{\infty,n}(t)\|_{L^2(\DD^n)} 
			\le C_{T,n,\alpha} e^{C_T nt}   k^{n/2+1} (C_T t^{(1-\theta/2)})^k (n+2k)^a	
\end{align}
			for all $k\ge 1$. Here $C_T$ and $a$ are constants independent of $n$ and $k$. We can fix a small time $t_*>0$ such that $C_T t_*^{(1-\theta/2)}<1$. Then for every fixed $t\le t_*$ and $n\in \mathbb{N}$, the right-hand side of \eqref{eq:final-Dyson} tends to zero as  $k\to \infty$, while the left-hand side is independent of $k$.  This yields $C_{\infty,n}(t) = 0$ for all $t\in [0,t_*]$ and for all $n\ge 1$. Repeating the above argument in several time intervals, we get $C_{\infty,n}(t) = 0$ for all time $t\in [0,T]$, yielding the uniqueness result. The proof of Theorem \ref{theo-main} is complete. 
\end{proof}

\end{document}